\documentclass[journal]{IEEEtran}
\ifCLASSINFOpdf

\else

\fi

\usepackage{mathtools}
\usepackage{psfrag}
\usepackage{rotating}
\usepackage[usenames]{color}
\usepackage{graphicx}
\graphicspath{{figures/}}
\usepackage{epsfig}    
\usepackage{amssymb,scalerel}   
\usepackage[font=footnotesize,caption=false]{subfig}
\usepackage{fixltx2e}
\usepackage{array}
\usepackage{dblfloatfix}

\usepackage{cite}

\usepackage{url}
\usepackage{tikz}

\usepackage{algorithm}     
\usepackage{algpseudocode} 

\DeclareMathOperator*{\minimize}{minimize}
\DeclareMathOperator*{\trace}{tr}

\DeclareMathOperator*{\minimum}{min}
\DeclareMathOperator*{\subject}{subject\ to}

\DeclareMathOperator*{\argmin}{arg\ min}

\DeclareMathOperator*{\diag}{diag}

\DeclareMathOperator*{\vectorize}{vec}

\DeclareMathOperator*{\parent}{par}

\DeclareMathOperator*{\children}{ch}

\DeclareMathOperator*{\svec}{svec}

\DeclareMathOperator*{\blkdiag}{blk\ diag}
\DeclareMathOperator*{\Ne}{Ne}

\DeclareMathOperator*{\leaves}{leaves}
\DeclareMathOperator*{\vectri}{vectri}

\newcounter{thm}

\newtheorem{theorem}[thm]{Theorem}

\newcounter{remcount}
\newtheorem{rem}[remcount]{Remark}

\definecolor{red}{rgb}{1,0,0}

\begin{document}
%
\title{Distributed Localization of Tree-structured Scattered Sensor Networks}
%
%
%
\author{Sina~Khoshfetrat~Pakazad,~\IEEEmembership{Member,~IEEE,}
        Emre~\"Ozkan,~\IEEEmembership{Member,~IEEE,}
        Carsten~Fritsche,~\IEEEmembership{Member,~IEEE,}
        Anders~Hansson,~\IEEEmembership{Member,~IEEE,}
        and~Fredrik~Gustafsson,~\IEEEmembership{Member,~IEEE}
\thanks{S. Khoshfetrat Pakazad, E. \"Ozkan, C. Fritsche, A. Hansson and F. Gustafsson are with the Division of Automatic Control, Department of Electrical Engineering, Link\"oping University, Sweden. Email: \{sina.kh.pa, emre, carsten, hansson, fredrik\}@isy.liu.se.}
}

\markboth{Journal of IEEE Transactions on Signal Processing}%
{Shell \MakeLowercase{\textit{et al.}}: Bare Demo of IEEEtran.cls for Journals}
%



\maketitle

\begin{abstract}
Many of the distributed localization algorithms are based on relaxed optimization formulations of the localization problem. These algorithms commonly rely on first-order optimization methods, and hence may require many iterations or communications among computational agents. Furthermore, some of these distributed algorithms put a considerable computational demand on the agents. In this paper, we show that for tree-structured scattered sensor networks, which are networks that their inter-sensor range measurement graphs have few edges (few range measurements among sensors) and can be represented using a tree, it is possible to devise an efficient distributed localization algorithm that solely relies on second-order methods. Particularly, we apply a state-of-the-art primal-dual interior-point method to a semidefinite relaxation of the maximum-likelihood formulation of the localization problem. We then show how it is possible to exploit the tree-structure in the network and use message-passing or dynamic programming over trees, to distribute computations among different computational agents. The resulting algorithm requires far fewer iterations and communications among agents to converge to an accurate estimate. Moreover, the number of required communications among agents, seems to be less sensitive and more robust to the number of sensors in the network, the number of available measurements and the quality of the measurements. This is in stark contrast to distributed algorithms that rely on first-order methods. We illustrate the performance of our algorithm using experiments based on simulated and real data.
\end{abstract}


%
\IEEEpeerreviewmaketitle

\section{Introduction}\label{sec:introduction}


The use of GPS for localizing sensor nodes in a sensor network is considered to be excessively expensive and wasteful, also in some cases intractable, \cite{bul:00,bis:06}. Instead many solutions for the localization problem tend to use inter-sensor distance or range measurements. In such a setting the localization problem is to find unknown locations of say $N$ sensors using existing noisy distance measurements among them and to sensors with known locations, also referred to as anchors. This problem is known to be NP hard \cite{mor:97}, and there have been many efforts to approximately solve this problem, \cite{kim:09,bis:04,wan:06,bis:06,gho:13,nar:14,cha:09,soa:15,sch:15,sri:08}.

One of the major approaches for approximating the localization problem, has been through the use of convex relaxation techniques, namely semidefinite, second-order and disk relaxations, see e.g., \cite{kim:09,bis:06,bis:04,wan:06,sri:08,gho:13,soa:15}. Although the centralized algorithms based on the these approximations reduce the computational complexity of solving the localization problem, they are still not scalable for solving large problems. Also centralized algorithms are generally communication intensive and more importantly lack robustness to failures. Furthermore, the use of these algorithms can become impractical due to certain structural constraints resulting from, e.g., privacy constraints and physical separation. These constraints generally prevent us from forming the localization problem in a centralized manner. One of the approaches to evade such issues is through the use of scalable and/or distributed algorithms for solving large localization problems. These algorithms enable us to solve the problem through collaboration and communication of several computational agents, which could correspond to sensors, without the need for a centralized computational unit. The design of distributed localization algorithms is commonly done by first reformulating the problem by exploiting or imposing structure on the problem and then employing efficient optimization algorithms for solving the reformulated problem, see e.g., some recent papers \cite{sim:14,gho:13,soa:15,sri:08}. For instance, authors in \cite{sri:08} put forth a solution for the localization problem based on minimization the discrepancy of the squared distances and the range measurements. They then propose a second-order cone relaxation for this problem and apply a Gauss-Seidel scheme to the resulting problem. This enables them to solve the problem distributedly. The proposed algorithm does not provide a guaranteed convergence and at each iteration of this algorithm, each agent is required to solve a second-order cone program, SOCP, which can potentially be expensive. Furthermore, due to the considered formulation of the localization problem, the resulting algorithm is prone to amplify the measurement errors and is sensitive to outliers. In \cite{sim:14}, the authors consider an SDP relaxation of the maximum likelihood formulation of the localization problem. They further relax the problem to an edge-based formulation as suggested in \cite{wan:06}. This then allows them to devise a distributed algorithm for solving the reformulated problem using alternating direction method of multipliers (ADMM). Even though this algorithm has convergence guarantees, each agent is required to solve an SDP at every iteration of the algorithm. In order to alleviate this, authors in \cite{gho:13} and \cite{soa:15} consider a disk relaxation of the localization problem and which correspond to an under-estimator of the original problem. They then use projection-based methods and Nestrov's optimal gradient method, respectively, for devising distributed algorithms for solving the resulting problem. These algorithms rely on finding a solution that lies in the intersection of the disks or spheres defined by the range measurements. Consequently, the computational demand on each agent for these algorithms is far less than the aforementioned algorithms. These algorithms commonly work well when there are many range measurements available and their performance is adversely affected if the number of measurements are decreased. Moreover, for the case of low quality, particularly biased, measurements, the convergence of the algorithms can be interrupted as the intersection can be empty.

The proposed algorithms in the aforementioned papers have been shown to be effective in analyzing large-scale localization problems. However, all these methods rely on first-order optimization algorithms and hence can require many iterations and communications to converge to an accurate enough solution. Furthermore, the number of iterations can vary significantly with different realizations of range measurements and changing topology of the sensor network. In this paper we show that in case it is possible to provide a tree representation of the inter-sensor range measurement graph of the sensor network (which is the case in many scenarios with few available range measurements), it is possible to alleviate these issues by devising far more efficient distributed localization algorithms that purely rely on second-order methods.

\subsection*{Contributions}

In this paper, we consider the localization problem for sensor networks where we have access to few range measurements among sensors. The availability of range measurements among $N$ sensors can be described using a graph with $N$ vertices or nodes and an edge between two nodes if there exists a range measurement between them. We refer to this graph as the inter-sensor measurement graph. For our purpose this graph is connected but sparse, i.e., it has few edges. For these sensors networks, it is commonly possible to represent the graph using a tree. We here propose a distributed localization algorithm based on the semidefinite relaxation of the localization problem \cite{sim:14,kim:09}. This algorithm relies on second-order methods, particularly state-of-the-art primal-dual interior-point methods, \cite{wri:97,tod:96,kho:15c,kho:15d}, and is obtained by distributing the computations of each iteration of the primal-dual method among several computational agents. This is done by first clustering the sensor nodes and providing a tree representation of inter-sensor measurement graph. The tree representation then allows us to use message-passing or dynamic programming over trees, \cite{kol:09,kho:15c,kho:15d,ber:73}, to compute the search directions at every iterations of the primal-dual methods distributedly by performing an upward-downward pass through the aforementioned tree. Consequently, and since primal-dual methods commonly converge within 20-50 iterations, our proposed algorithm in comparison to existing ones requires far fewer iterations and communications among agents to converge to a solution. Furthermore, the computational burden for each agent at each iteration only concerns factorizing a relatively small matrix, c.f., \cite{sim:14,sri:08}.

\subsection*{Outline}

In Section \ref{sec:MLL} we review a maximum-likelihood formulation of the localization problem. Section \ref{sec:scattered} provides a formal description of tree-structured scattered sensor networks and describes how the structure in the problem can be reflected in the localization optimization problem. Section \ref{sec:chordal} reviews how certain structure in nonlinear SDPs enable us to utilize domain-space decomposition to decompose them. This decomposition technique is then used in Section \ref{sec:decompositionloc} to decompose the localization optimization problem. In this section we also describe how the decomposed problem can be written as coupled SDP. We then put forth a generic description of primal-dual interior-point methods in Section \ref{sec:DPDIPM} and show how they, in combination with message-passing, can be used to devise efficient distributed solvers for the localization problems. In this section we also discuss the computational and communication complexity of the proposed distributed algorithm. The numerical experiments are presented in Section~\ref{sec:numerical}, and we conclude the paper with final remarks in Section~\ref{sec:conclusions}.

\subsection*{Notations and Definitions}\label{sec:notation}
We denote by $\mathbb R$ the set of real scalars and by $\mathbb R^{n\times m}$ the set of real $n\times m$ matrices. The set of $n \times n$ symmetric matrices are represented by $\mathbf S^n$. The transpose of a matrix $A$ is denoted by $A^T$ and the column and null space of this matrix is denoted by $\mathcal{C}(A)$ and $\mathcal N(A)$, respectively. We denote the set of positive integers $\{1,2,\ldots,p\}$ with $\mathbb{N}_p$. Given a set $J \subset \mathbb{N}_n$, the matrix $E_J \in \mathbb{R}^{|J|\times n}$ is the $0$-$1$ matrix that is obtained by deleting the rows indexed by $\mathbb{N}_n \setminus J$ from an identity matrix of order $n$, where $|J|$ denotes the number of elements in set $J$. This means that $E_Jx$ is a $|J|$- dimensional vector with the components of $x$ that correspond to the elements in $J$, and we denote this vector with $x_J$. Also $e_j$ denotes a 0--1 $n$-dimensional  vector with only a nonzero element at the $j$th component. Similarly, given $J \subset \mathbb N_n$, $e_J$ denotes a 0--1 $n$-dimensional vector with ones at elements specified by $J$. With $x^{i,(k)}_l$ we denote the $l$th element of vector $x^i$ at the $k$th iteration. Also given vectors $x^i$ for $i= 1, \dots, N$, the column vector $(x^1, \dots, x^N)$ is all of the given vectors stacked. For a vector $x$, with $\diag(x)$ we denote a diagonal matrix with its diagonal elements given by $x$. Similarly, given matrices $X^i$ for $i = 1, \dots, N$, with $\blkdiag(X^1, \dots, X^N)$ we denote a block-diagonal matrix with diagonal blocks given by each of the given matrices. For a matrix $X \in \mathbb R^{m\times n}$, $\vectorize(X)$ is an $mn$-dimensional vector that is obtained by stacking all columns of $X$ on top of each other. Given a symmetric matrix $X \in \mathbf S^n$
\begin{multline*}
\svec(X) := (X_{11}, \sqrt{2}X_{21}, \dots, \sqrt{2} X_{n1}, X_{22},\\ \sqrt{2} X_{32}, \dots, \sqrt{2} X_{n2}, \dots, X_{nn}).
\end{multline*}
Also for a square matrix $X \in \mathbb R^{n\times n}$ we denote with $\vectri(X)$ a column vector which includes all elements on the upper triangle of $X$ stacked. Given two matrices $X$ and $Y$ by $X\otimes Y$ we denote the standard Kronecker product. Given $X \in \mathbf S^n$, define $U$ as an $n(n + 1)/2 \times n^2$ matrix such that $U \vectorize(X) = \svec(X)$. Then for two matrices $X, Y \in \mathbb R^{n\times n}$, $\otimes_s$ denotes the symmetrized Kronecker product that is defined as
\begin{align*}
X \otimes_s Y := \frac{1}{2} U(X \otimes Y + Y \otimes X)U^T.
\end{align*}
For properties of the symmetrized Kronecker product refer to \cite{tod:96}.

A graph is denoted by $Q(V,\mathcal E)$ where $V = \{1, \dots, n\}$ is its set of vertices or nodes and $\mathcal E \subseteq V\times V$ denotes its set of edges. Vertices $i, j \in V$ are adjacent if $(i, j) \in E$, and we denote the set of adjacent vertices of $i$ by $\Ne(i) = \{ j \in V | (i, j) \in\mathcal E \}$. A graph is said to be complete if all its vertices are adjacent. An induced graph by $V^\prime \subseteq V$ on $Q(V,\mathcal E)$, is a graph $Q_I(V^\prime,\mathcal E^\prime)$ where $\mathcal E^\prime = \mathcal E\cap V^\prime \times V^\prime$. A clique $C_i$ of $Q(V,\mathcal E)$ is a maximal subset of $V$ that induces a complete subgraph on $Q$, i.e., no clique is properly contained in another clique, \cite{blp:94}. Assume that all cycles of length at least four of $Q(V,\mathcal E)$ have a chord, where a chord is an edge between two non-consecutive vertices in a cycle. This graph is then called chordal \cite[Ch. 4]{gol:04}. It is possible to make graphs chordal by adding edges to the graph. The resulting graph is then referred to as a chordal embedding. Let $\mathbf C_Q = \{ C_1, \dots, C_q \}$ denote the set of its cliques, where $q$ is the number of cliques of the graph. Then there exists a tree defined on $\mathbf C_Q$ such that for every $C_i, C_j \in\mathbf C_Q$ where $i \neq j$, $C_i \cap C_j$ is contained in all the cliques in the path connecting the two cliques in the tree. This property is called the clique intersection property, \cite{blp:94}. Trees with this property are referred to as clique trees.
\section{Maximum Likelihood Localization}\label{sec:MLL}
In this paper we consider a localization problem for a network of $N$ sensors distributed in an area in presence of $m$ anchors. The exact locations of these sensors, $x^i_s$, are deemed to be unknown however we assume that the positions of the anchors, $x^i_a$, are given. Furthermore, the sensors are capable of performing computations and some can measure their distance to certain sensors and some of the anchors. We assume that if sensor $i$ can measure its distance to sensor $j$ so can sensor $j$ measure its distance to sensor $i$. This then allows us to describe the range measurement availability among sensors using an undirected graph $G_r(V_r, \mathcal E_r)$ with vertex set $V_r = \{ 1, \dots, N \}$ and edge set $\mathcal E_r$. An edge $(i,j) \in \mathcal E_r$ if and only if a range measurement between sensors $i$ and $j$ is available. We refer to this graph as inter-sensor measurement graph and assume that it is connected. Let us define the set of neighbors of each sensor $i$, $ \text{Ne}_r(i)$, as the set of sensors to which this sensor has an available range measurement. In a similar fashion let us denote the set of anchors to which sensor $i$ can measure its distance to by $\text{Ne}_a(i) \subseteq \{ 1, \dots, m\}$. Let us describe the inter-sensor range measurements for each sensor, $i \in \mathbb N_N$, as
\begin{align}\label{eq:InterSensorMeasurement}
\mathcal R_{ij} = \mathcal D_{ij} + E_{ij}, \quad  j \in \text{Ne}_r(i),
\end{align}
where $\mathcal D_{ij} = \| x_s^i - x_s^j \|_2 $ defines the noise-free sensor distance, $E_{ij}$ is the inter-sensor measurement noise and $E_{ij} \sim P_{ij}^s(\mathcal D_{ij}|\mathcal R_{ij})$ with $P_{ij}^s(\cdot)$ being the so-called inter-sensor sensing probability density function (PDF). We here make the standard assumption that $\mathcal R_{ij} = \mathcal R_{ji}$, see e.g., \cite{soa:15,shi:10}. Similarly we can describe the anchor range measurements for each sensor $i$ as
\begin{align}\label{eq:AnchorSensorMeasurement}
\mathcal Y_{ij} = \mathcal Z_{ij} + V_{ij}, \quad  j \in \text{Ne}_a(i),
\end{align}
where $\mathcal Z_{ij} = \| x_s^i - x_a^j \|_2$ defines the noise-free anchor-sensor distance, $V_{ij}$ is the anchor-sensor measurement noise and $V_{ij} \sim P_{ij}^a(\mathcal Z_{ij}|\mathcal Y_{ij})$ with $P_{ij}^a(\cdot)$ being the so-called anchor-sensor sensing PDF. Here we assume that the inter-sensor and anchor-sensor measurement noise PDFs, i.e., $P_{ij}^s(\cdot)$ and $P_{ij}^a(\cdot)$, respectively, are Gaussian. Particularly, we assume that the inter-sensor and anchor-sensor measurement noises are independent and that $E_{ij} \sim \mathcal N(0, \Sigma^r_{ij})$ and $V_{ij} \sim \mathcal N(0, \Sigma^a_{ij})$. Notice that this assumption can be relaxed to any distribution that is a log-concave function of distances $\mathcal D_{ij}$ and $\mathcal Z_{ij}$, however, for the sake of brevity we limit ourselves to the case of Gaussian distributions. Having defined the setup of the sensor network, we can write the localization problem in a maximum likelihood setting as

\small
\begin{multline}\label{eq:MLOriginal}
X_{\text{ML}}^\ast = \argmin_{X} \Bigg\{ \sum_{i = 1}^N\Bigg( \sum_{\tiny\begin{split}j \in &\text{Ne}_r(i)\\ i &< j\end{split}\normalsize} \frac{1}{\Sigma^r_{ij}}\left( \mathcal D_{ij}(x^i_s, x^j_s) - \mathcal R_{ij}\right)^2  \\  + \sum_{j \in \text{Ne}_a(i)} \frac{1}{\Sigma^a_{ij}}\left(\mathcal Z_{ij}(x^i_s, x^j_a) - \mathcal Y_{ij}\right)^2 \Bigg) \Bigg\},
\end{multline}
\normalsize
where $X = \begin{bmatrix} x^1_s & \dots & x^N_s \end{bmatrix}\in \mathbb R^{d\times N}$ with $d = 2$ or $d = 3$. This problem can be formulated as a constrained optimization problem, as was described in \cite{sim:14}, which is discussed next.
First let us define the function

\small
\begin{multline}
f(\Lambda, \Xi, D, Z) : =  \sum_{i = 1}^N  \Bigg(\sum_{\tiny\begin{split}j \in &\text{Ne}_r(i)\\ i &< j\end{split}\normalsize}\frac{1}{\Sigma^r_{ij}} (\Lambda_{ij} - 2D_{ij}\mathcal R_{ij} +\mathcal R_{ij}^2) \\  + \sum_{j \in \text{Ne}_a(i)}\frac{1}{\Sigma^a_{ij}} (\Xi_{ij} - 2Z_{ij}\mathcal Y_{ij} +\mathcal Y_{ij}^2)  \Bigg).
\end{multline}
\normalsize
Then the problem in \eqref{eq:MLOriginal} can be equivalently rewritten as the following constrained optimization problem

\small
\begin{subequations}\label{eq:MLConstrained}
\begin{align}
\minimize_{X, S, \Lambda, \Xi, D, Z}& \quad f(\Lambda, \Xi, D, Z)\\
\subject & \notag \\ & \begin{rcases*} S_{ii} + S_{jj} -2S_{ij} = \Lambda_{ij} \\ \Lambda_{ij} = D_{ij}^2, \quad D_{ij} \geq 0, \ j\in \text{Ne}_r(i),i<j\end{rcases*}, \ i \in \mathbb N_N\label{eq:MLConstrained-b}\\ &\begin{rcases*}  S_{ii} - 2(x_s^i)^Tx_a^j + \| x_a^j \|^2_2  = \Xi_{ij} \\ \Xi_{ij} = Z_{ij}^2, \quad Z_{ij} \geq 0, \ \ j \in \text{Ne}_a(i) \end{rcases*}, \ \ i \in \mathbb N_N \label{eq:MLConstrained-c}\\
&  S = X^T X.\label{eq:MLConstrained-d}
\end{align}
\end{subequations}
\normalsize
So far we have reviewed a way to formulate the localization problem over general sensor networks as a constrained optimization problem. In this paper, however, we are particularly interested in localization problem pertaining to sensor networks with an inherent tree structure which relies on the assumption that the graph $G_r(V_r, \mathcal E_r)$ can be represented using a tree. We describe the localization problem of such networks in the next section.

\section{Localization of Tree-structured Scattered Sensor Networks}\label{sec:scattered}

Let the graph $G_r(V_r, \mathcal E_r)$ be connected with few edges. Also assume that a chordal embedding $\bar G_r(V_r, \bar{\mathcal E}_r)$ of this graph can be achieved by adding only a few edges. This graph can then be represented using its clique tree. Furthermore, given the set of its cliques $\mathbf C_{\bar G_r} = \{ C_1, \dots, C_q \}$, we have $| C_i | \ll N$. We refer to such sensor networks as tree-structured scattered. The localization problem of these sensor networks can also be formulated as a constrained optimization problem using the approach discussed in Section \ref{sec:notation}. However, the formulation of the problem in~\eqref{eq:MLConstrained} is not fully representative of the structure in the problem. In order to exploit the structure in our localization problem we modify \eqref{eq:MLConstrained}, and equivalently rewrite it as

\small
\begin{subequations}\label{eq:MLConstrainedScattered}
\begin{align}
\minimize_{X, S, \Lambda, \Xi, D, Z}& \quad f(\Lambda, \Xi, D, Z) \\
\subject & \notag \\ & \  \begin{rcases*} S_{ii} + S_{jj} -2S_{ij} = \Lambda_{ij} \\ \Lambda_{ij} = D_{ij}^2, \quad D_{ij} \geq 0,  \ j\in \text{Ne}_r(i), i<j\end{rcases*}, \ i \in \mathbb N_N\label{eq:MLConstrainedScattered-b}\\ &\  \begin{rcases*}  S_{ii} - 2(x_s^i)^Tx_a^j + \| x_a^j \|^2_2  = \Xi_{ij} \\ \Xi_{ij} = Z_{ij}^2, \quad Z_{ij} \geq 0, \ \ j \in \text{Ne}_a(i) \end{rcases*},  \ i \in \mathbb N_N \label{eq:MLConstrainedScattered-c}\\
& \   S \succeq 0, \ \  S_{ij} = (x_s^i)^Tx_s^j, \notag \\ &\hspace{25mm} \forall \ (i,j) \in \mathcal E_r \cup \{ (i, i) \ | \ i \in V_r \}\label{eq:MLConstrainedScattered-d}.
\end{align}
\end{subequations}
\normalsize
Note that, here, we have modified the constraint in \eqref{eq:MLConstrained-d} so that the structure in the problem is more explicit. This modification is based on the observation that not all the elements of $S$ are used in \eqref{eq:MLConstrained-b} and \eqref{eq:MLConstrained-c}, and hence we only have to specify the ones that are needed and can leave the rest free. In \cite{kim:09}, \cite{wan:06}, the authors first conduct a semidefinite relaxation on \eqref{eq:MLConstrained}. They then exploit the structure as we did in \eqref{eq:MLConstrainedScattered} and use the ideas in \cite{fukuda_exploitingsparsity} to devise efficient centralized solvers for the localization problem. Here, however, we stick to the formulation in \eqref{eq:MLConstrainedScattered} which is a nonlinear SDP, and use scheme in \cite{kim+koj+mev+yam10} to decompose this problem directly. We then perform a semidefinite relaxation on the resulting problem and rewrite the problem as a coupled SDP. This in turn facilitates the use of efficient scalable or distributed solvers. The use of the so-called domain-space decomposition presented in \cite{kim+koj+mev+yam10} is at the heart of this reformulation approach. We review this decomposition scheme next, for the sake of completeness.
\begin{rem}
Notice that the added edges for computing a chordal embedding for the inter-sensor measurement graph does not affect the problem description in \eqref{eq:MLConstrainedScattered}, and only facilitates the clustering of the sensor nodes.
\end{rem}
\section{Chordal Sparsity in Semidefinite Programs}\label{sec:chordal}
In this section we first briefly review some of important properties of sparse semidefinite matrices and then discuss how these can be used for reformulating semidefinite programs with chordal sparsity suitable to be solved distributedly.

\subsection{Chordal Sparsity}\label{sec:sparsity}

Graphs can be used to characterize partial symmetric matrices. Partial symmetric matrices correspond to symmetric matrices where only a subset of their elements are specified and the rest are free. We denote the set of all $n \times n$ partially symmetric matrices on a graph $Q(V,\mathcal E)$ by $\mathbf S_Q^n$, where only elements with indices belonging to $\mathbf I_s = \mathcal E \cup \{ (i,i) \ | \ i \in \mathbb N_n\}$ are specified. Now consider a matrix $X \in \mathbf S_Q^n$. Then $X$ is positive semidefinite completable if by manipulating its free elements, i.e., elements with indices belonging to $\mathbf I_f = (V \times V) \setminus \mathbf I_s $, we can generate a positive semidefinite matrix. The following theorem states a fundamental result on positive semidefinite completion.
\begin{theorem}(\cite[Thm. 7]{gro:84})\label{thm:NSDC}
Let $Q(V,\mathcal E)$ be a chordal graph with cliques $C_1, \dots, C_q$ such that clique intersection property holds. Then $X \in \mathbf S_Q^n$ is positive semidefinite completable, if and only if
\begin{align}
X_{C_iC_i} \succeq 0, \quad   \ i \in \mathbb N_q,
\end{align}
where $X_{C_iC_i} = E_{C_i} X E_{C_i}^T$.
\end{theorem}

Note that the matrices $X_{C_iC_i}$ for $i \in \mathbb N_q$, are the fully specified principle submatrices of $X$. Hence, Theorem~\ref{thm:NSDC} states that a chordal matrix $X \in \mathbf S_{Q}^n$ is positive  semidefinite completable if and only if all its fully specified principle submatices are positive semidefinite. As we will see next this property can be used for decomposing SDPs with this structure.

\subsection{Domain-space Decomposition}\label{sec:decomposition}

Consider a chordal graph $Q(V,\mathcal E)$, with $\{ C_1, \dots, C_q \}$ the set of cliques such that the clique intersection property holds. Let us define sets $J_i \subset \mathbb N_n$ such that the sparsity pattern graph for $\sum_{i=1}^N e_{J_i}e_{J_i}^T $ is $Q(V,\mathcal E)$. Then for the following nonlinear SDP
\begin{subequations}\label{eq:SDP}
\begin{align}
 \minimize_{z^1, \dots, z^N, X} &\quad  \sum_{i=1}^N  f^i(z^i,\svec(E_{J_i}XE_{J_i}^T)) \label{eq:SDP1} \\
 \subject & \quad g^i(z^i, \svec(E_{ J_i}XE_{J_i}^T)) \in \Omega_i, \hspace{3mm} i \in \mathbb N_N, \label{eq:SDP2} \\
& \quad X \succeq 0, \label{eq:SDP3}
\end{align}
\end{subequations}
the only elements of $X$ that affect the cost function in \eqref{eq:SDP1} and the constraint in \eqref{eq:SDP2} are elements specified by indices in $\mathbf I_s$. Using Theorem~\ref{thm:NSDC}, the optimization problem in~\eqref{eq:SDP} can then be equivalently rewritten as
\begin{subequations}\label{eq:MSDP1}
\begin{align}
 \minimize_{z^1, \dots, z^N, X} &  \quad \sum_{i=1}^N  f^i(z^i,\svec(E_{J_i}XE_{J_i}^T)) \label{eq:MSDP1-1} \\
 \subject  & \quad g^i(z^i, \svec(E_{ J_i}XE_{J_i}^T)), \hspace{3mm} i \in \mathbb N_N,  \label{eq:MSDP1-2} \\
& \quad  X_{C_iC_i} \succeq 0, \hspace{3mm} i \in \mathbb N_q,  \label{eq:MSDP1-3}
\end{align}
\end{subequations}
where notice that the constraints in \eqref{eq:MSDP1-3} are coupled semidefinite constraints, \cite{fukuda_exploitingsparsity, kim+koj+mev+yam10}. It is possible to explicitly describe the coupling using consistency constraints and rewrite \eqref{eq:MSDP1} as
\begin{subequations}\label{eq:MSDP}
\begin{align}
 \minimize_{z^1, \dots, z^N, X^1, \dots, X^q, X} &  \quad \sum_{i=1}^N  f^i(z^i,\svec(E_{J_i}XE_{J_i}^T)) \label{eq:MSDP-1} \\
 \subject  & \quad g^i(z^i, \svec(E_{ J_i}XE_{ J_i}^T)), \hspace{3mm} i \in \mathbb N_N, \label{eq:MSDP-2} \\
& \quad  X^{i} \succeq 0, \hspace{3mm} i \in \mathbb N_q,  \label{eq:MSDP-3}\\
& \quad  X^i = E_{C_i}XE_{C_i}^T, \quad i \in \mathbb N_q,  \label{eq:MSDP-5}
\end{align}
\end{subequations}
where $X^i \in \mathbf S^{|C_i|}$. This method of reformulating \eqref{eq:SDP} as \eqref{eq:MSDP} is referred to as the domain-space decomposition, \cite{kim+koj+mev+yam10,and:11}. The structure in the localization of tree-structured scattered sensor networks enable us to use this technique for reformulating the problem in such a way that would better facilitate the use of efficient distributed solvers. This is discussed in the next section.

\section{Decomposition and Convex Formulation of Localization of Tree-structured Scattered Sensor Networks}\label{sec:decompositionloc}

Consider the inter-sensor measurement graph $G_r(V_r, \mathcal E_r)$, and assume that it is chordal. In case this graph is not chordal the upcoming discussions hold for any of its chordal embeddings. Let $\mathbf C_{G_r} = \{ C_1, \dots, C_q \}$ and $T(V_t,\mathcal E_t)$ be a clique tree. Based on the discussion in Section~\ref{sec:decomposition}, then for the problem in \eqref{eq:MLConstrainedScattered} we have $S \in \mathbf S^N_{G_r}$. Hence, we can rewrite~\eqref{eq:MLConstrainedScattered} as

\small
\begin{subequations}\label{eq:MLConstrainedScattered1}
\begin{align}
\minimize_{X, S_{C_iC_i}, \Lambda, \Xi, D, Z}& \quad f(\Lambda, \Xi, D, Z)\label{eq:MLConstrainedScattered1-a} \\
\subject & \notag \\ & \begin{rcases*} S_{ii} + S_{jj} -2S_{ij} = \Lambda_{ij} \\ \Lambda_{ij} = D_{ij}^2, \quad D_{ij} \geq 0, j\in \text{Ne}_r(i), i<j\end{rcases*}, i \in \mathbb N_N,\label{eq:MLConstrainedScattered1-b}\\ & \begin{rcases*}  S_{ii} - 2(x_s^i)^Tx_a^j + \| x_a^j \|^2_2  = \Xi_{ij} \\ \Xi_{ij} = Z_{ij}^2, \quad Z_{ij} \geq 0, \ \ j \in \text{Ne}_a(i) \end{rcases*}, \ \ i \in \mathbb N_N, \label{eq:MLConstrainedScattered1-c}\\
&  S_{C_iC_i} \succeq 0, \quad S_{C_iC_i} = E_{C_i} X^TXE_{C_i}^T,   \ \  i \in \mathbb N_q  \label{eq:MLConstrainedScattered1-d},
\end{align}
\end{subequations}
\normalsize
Notice that even though the cost function for this problem is convex, the constraints in \eqref{eq:MLConstrainedScattered1-b}--\eqref{eq:MLConstrainedScattered1-d} are non-convex and hence the problem is non-convex. Consequently, we next address the localization problem by considering a convex relaxation of this problem. This allows us to solve the localization problem approximately.

One of the ways to provide a convex approximation of the problem in \eqref{eq:MLConstrainedScattered1} is to relax the quadratic equality constraints in \eqref{eq:MLConstrainedScattered1-b}--\eqref{eq:MLConstrainedScattered1-d} using Schur complements, which results in

\small
\begin{subequations}\label{eq:MLConstrainedScattered2}
\begin{align}
&\minimize_{\tiny\begin{matrix}X, S_{C_iC_i}, \Lambda_{ij}, \Xi_{ij}, D_{ij},\\ Z_{ij}, T^i, \Gamma^{ij}, \Phi^{ij} \end{matrix}\normalsize} \quad \sum_{i = 1}^N  \left(\sum_{\tiny\begin{matrix}j \in \text{Ne}_r(i)\\ i < j\end{matrix}\normalsize}f_{ij}(\Lambda_{ij}, D_{ij}) +\right. \notag\\& \left. \hspace{50mm} \sum_{j \in \text{Ne}_a(i)}g_{ij}(\Xi_{ij}, Z_{ij})  \right) \label{eq:MLConstrainedScattered2-a} \\
&\subject \notag\\ & \quad\quad \quad(S_{ii}, S_{jj}, S_{ij}, \Lambda_{ij}, D_{ij}, \Gamma^{ij}) \in \Omega_{ij}, \ (i,j) \in \mathcal E_r,\ i<j, \label{eq:MLConstrainedScattered2-b}\\
&\quad\quad\quad (S_{ii}, x_s^i, \Xi_{ij} , Z_{ij},\Phi^{ij}) \in \Theta_{ij}, \ \ j \in \text{Ne}_a(i), \ \ i \in \mathbb N_N, \label{eq:MLConstrainedScattered2-c}\\
& \quad\quad\quad \begin{bmatrix} I & XE_{C_i}^T \\ E_{C_i}X^T & S_{C_iC_i} \end{bmatrix} = T^i, \quad T^i \succeq 0, \ \  i \in \mathbb N_q  \label{eq:MLConstrainedScattered2-d},
\end{align}
\end{subequations}
\normalsize
where

\small
\begin{align*}
f_{ij}(\Lambda_{ij}, D_{ij}) & = \frac{1}{\sigma_{ij}^2} (\Lambda_{ij} - 2D_{ij}R_{ij} + R_{ij}^2),\\
g_{ij}(\Xi_{ij}, Z_{ij}) & = \frac{1}{\delta_{ij}^2} (\Xi_{ij} - 2Z_{ij}Y_{ij} + Y_{ij}^2),
\end{align*}
\normalsize
and

\small
\begin{align*}
\Omega_{ij} &= \Bigg \{ (S_{ii}, S_{jj}, S_{ij}, \Lambda_{ij}, D_{ij}, \Gamma^{ij}) \ \Bigg | \ S_{ii} + S_{jj} -2S_{ij} = \Lambda_{ij},  \\ & \quad \quad \quad \quad \quad \quad \quad \quad \quad    \ \begin{bmatrix} 1 & D_{ij} \\ D_{ij} & \Lambda_{ij} \end{bmatrix} = \Gamma^{ij}, \ \Gamma^{ij} \succeq 0, \ D_{ij} \geq 0 \Bigg\},\\
\Theta_{ij} &= \Bigg \{ (S_{ii}, x_s^i, \Xi_{ij} , Z_{ij}, \Phi^{ij}) \ \Bigg | \ S_{ii} - 2(x_s^i)^Tx_a^j + \| x_a^j \|^2_2  = \Xi_{ij},  \\ & \quad \quad \quad   \ \  \ \begin{bmatrix} 1 & Z_{ij} \\ Z_{ij} & \Xi_{ij} \end{bmatrix} = \Phi^{ij}, \ \Phi^{ij} \succeq 0 , Z_{ij} \geq 0, \ \ j \in \text{Ne}_a(i) \Bigg\},
\end{align*}
\normalsize
with the variables $\Gamma^{ij}$, $\Phi^{ij}$ and $T^i$ as slack variables. The addition of the slack variables enable us to make the description of the semidefinite constraints simpler. This problem is a coupled SDP and can be solved distributedly using $q$ computational agents. In order to see this with more ease, let us introduce a grouping of the cost function terms and constraints in~\eqref{eq:MLConstrainedScattered2-a}--\eqref{eq:MLConstrainedScattered2-c}. To this end we first describe a set of assignment rules. It is possible to assign
\begin{enumerate}
\item the constraint $(S_{ii}, S_{jj}, S_{ij}, \Lambda_{ij}, D_{ij}, \Gamma^{ij}) \in \Omega_{ij}$ and the cost function term $f_{ij}$ to agent $k$ if $(i,j) \in C_k \times C_k$;
\item the set of constraints $(S_{ii}, x_s^i, \Xi_{ij} , Z_{ij}, \Phi^{ij}) \in \Theta_{ij}, \ \ j \in \text{Ne}_a(i)$ and the cost function terms  $g_{ij}, \ \ j \in \text{Ne}_a(i)$ to agent $k$ if $i \in C_k$.
\end{enumerate}
We denote the indices of the constraints and cost function terms assigned to agent $k$ through Rule 1 above as $\phi_{k}$, and similarly we denote the set of constraints and cost function terms that are assigned to agent $k$ through Rule 2 by $\bar \phi_k$. Using the mentioned rules and the defined notations, we can now group the constraints and the cost function terms and rewrite the problem in \eqref{eq:MLConstrainedScattered2} as

\small
\begin{subequations}\label{eq:MLConstrainedScattered3}
\begin{align}
&\minimize_{\tiny\begin{matrix}X, S_{C_iC_i}, \Lambda_{ij}, \Xi_{ij}, D_{ij},\\ Z_{ij}, T^i, \Gamma^{ij}, \Phi^{ij} \end{matrix}\normalsize} \quad \sum_{k = 1}^q \left[ \sum_{(i,j) \in \phi_k} f_{ij}(\Lambda_{ij}, D_{ij}) +  \right. \notag \\ &\left. \hspace{50mm} \sum_{i \in \bar \phi_k} \sum_{j \in \text{Ne}_a(i)} g_{ij}(\Xi_{ij}, Z_{ij}) \right] \label{eq:MLConstrainedScattered3-a} \\
&\subject  \notag \\ & \quad \quad  \begin{rcases*}   (S_{ii}, S_{jj}, S_{ij}, \Lambda_{ij}, D_{ij}, \Gamma^{ij}) \in \Omega_{ij}, \ \ (i,j) \in \phi_k \\  (S_{ii}, x_s^i, \Xi_{ij} , Z_{ij},\Phi^{ij}) \in \Theta_{ij}, \ \ j \in \text{Ne}_a(i) \ \ i \in \bar \phi_k \\ \begin{bmatrix} I & XE_{C_k}^T \\ E_{C_k}X^T & S_{C_kC_k} \end{bmatrix} = T^k, \ T^k \succeq 0 \end{rcases*},  k \in \mathbb N_q  \label{eq:MLConstrainedScattered3-b}
\end{align}
\end{subequations}
\normalsize
Notice that this problem can now be seen as a combination of $q$ coupled subproblems, each defined by a term in the cost function together with its corresponding set of constraints in~\eqref{eq:MLConstrainedScattered3-b}. It is possible to decompose this problem by introducing additional local variables and consistency constraints and use any proximal point splitting method, e.g., ADMM, to solve this problem distributedly. However, there are major disadvantages for the resulting distributed solution, such as
\begin{itemize}
\item the local subproblems that needs to be solved by each agent is a semidefinite program that are computationally expensive to solve;
\item inexact solutions for semidefinite programs can be far away from the optimal solution;
\item the algorithm generally requires many iterations to converge to an accurate solution that particularly satisfies the consistency constraints;
\item the number of consistency constraints are generally big for such problems which can even further adversely affect the convergence and numerical properties of such algorithms.
\end{itemize}
In order to evade the aforementioned issues, we next put forth an alternative distributed algorithm based on primal-dual interior-point methods that fully takes advantage of the structure in the problem and yields an accurate solution within much lower number of iterations and with far less computational demands from each agent.

\begin{rem}
The accuracy of the estimates obtained from solving \eqref{eq:MLConstrainedScattered3} can be improved by pushing the rank of matrices $\Gamma^{ij}$ and $\Phi^{ij}$ to 1 and the rank of matrices $T^i$ to $d$, see e.g., \cite{wan:06}. One way to achieve this is through the use of nuclear norm regularization by adding
\small
\begin{align}
\sum_{k = 1}^q \left[ \alpha^k \| T^k \|_* + \sum_{(i,j) \in \phi_k} \rho^{ij} \| \Gamma^{ij} \|_* + \sum_{i \in \bar \phi_k} \sum_{j \in \text{Ne}_a(i)} \mu^{ij} \| \Phi^{ij} \|_* \right],
\end{align}
\normalsize
to the cost function of \eqref{eq:MLConstrainedScattered3}, see \cite{rec:10}, where $\| \cdot \|_*$ denotes the nuclear norm of a matrix and $\alpha^k>0$, $\rho^{ij}>0$ and $\mu^{ij}>0$ are the so-called regularization parameters. Since all the aforementioned matrices are restricted to be positive semidefinite this will be equivalent to
\small
\begin{align}
\sum_{k = 1}^q \left[ \alpha^k \trace( T^k ) + \sum_{(i,j) \in \phi_k} \rho^{ij} \trace(  \Gamma^{ij} ) + \sum_{i \in \bar \phi_k} \sum_{j \in \text{Ne}_a(i)} \mu^{ij} \trace(  \Phi^{ij} ) \right].
\end{align}
\normalsize
Notice that by increasing the regularization parameters the rank of these matrices are further pushed towards lower values. Furthermore, this does not affect the coupling structure in the problem since the added terms to the cost function concern the local matrix variables. Here, for the sake of brevity and notational simplicity, we do not consider the use of regularization. The coming discussion in Section \ref{sec:DPDIPM} can be extended to the regularized problem with little effort.
\end{rem}

\subsection{A Simple Assignment Strategy}

Before we continue, let us first put forth an assignment strategy that is simple and satisfies the assignment rules discussed above. Recall that in order to form the problem in \eqref{eq:MLConstrainedScattered3}, we first need to cluster the sensor nodes. Based on this clustering, we use the assignment strategy described in Algorithm \ref{alg:Ass}.
\begin{algorithm}[tb]
\caption{A Simple Assignment Strategy}\label{alg:Ass}
\begin{algorithmic}[1]
\small
\State{Given the inter-sensor measurement graph $G_r(V_r \mathcal E_r)$ and $C_{G_r} = \{ C_1, \dots, C_q \}$}
\For {$k = 1, \dots, q$}
\For {$i \in C_k$}
\For {$j \in \text{Ne}_r(i)$ and $i<j$}
\If{$\Omega_{ij}$ is not assigned and $j\in C_k$}
\State{Assign it to agent $k$}
\EndIf
\If{$f_{ij}$ is not assigned and $j\in C_k$}
\State{Assign it to agent $k$}
\EndIf
\EndFor
\For {$j \in \text{Ne}_a(i)$}
\If{$\Theta_{ij}$ is not assigned}
\State{Assign it to agent $k$}
\EndIf
\If{$g_{ij}$ is not assigned}
\State{Assign it to agent $k$}
\EndIf
\EndFor
\EndFor
\EndFor
\normalsize
\end{algorithmic}
\end{algorithm}
Notice that the resulting assignment heavily relies on the ordering of the cliques or clusters of sensors. Consequently, different ordering of the cliques may result in different assignments of constraints and terms in the objective function. Furthermore, even though this assignment algorithm is simple, it may lead to unbalanced distribution of constraints and cost function terms. This means that some agents maybe assigned a disproportionate number of variables, constraints and objective function terms. One can avoid such a situation by modifying the \emph{if} statements in steps 5, 8, 13 and 16 of the algorithm, by adding watchdogs that prevent unbalanced assignments. For the sake of brevity and so as to not clutter the presentation, we do not discuss this any further.
\begin{rem}
Notice that each pair $\Omega_{ij}$ and $f_{ij}$ corresponds to the range measurement between sensors $i$ and $j$ and each pair $\Theta_{ij}$ and $g_{ij}$ corresponds to a range measurement between sensor $i$ and anchor $j$. Based on this, using the assignment rules, we essentially assign different range measurements to each sensor cluster or computational agent.
\end{rem}
\section{Distributed Primal-dual Interior-point Method for Coupled SDPs}\label{sec:DPDIPM}
The problem in \eqref{eq:MLConstrainedScattered3} can be written in the following standard form

\small
\begin{subequations}
\begin{align}
\minimize &\quad \sum_{i = 1}^q (c^i)^T y \\
\subject &\notag \\ & \begin{rcases*} Q^i_j \svec(X_j^i) + W_j^i y = b^i_j, \quad j = 1, \dots, m_i \\ A^iy = \bar b^i \\ D^iy \leq g^i \\ X^i_j \succeq 0, \quad j = 1, \dots, m_i \end{rcases*}, \ i \in \mathbb N_q
\end{align}
\end{subequations}
\normalsize
where the variables $X^i_j$ and $y$ are matrix and linear variables, respectively. This problem can be written more compactly as
\begin{subequations}\label{eq:GeneralSDP}
\begin{align}
\minimize &\quad \sum_{i = 1}^q (c^i)^T y \\
\subject &\notag \\ & \quad \begin{rcases*} Q^i x^i + W^i y = b^i \\ A^iy = \bar b^i \\ D^iy \leq g^i \\ X^i \succeq 0 \end{rcases*}, \ i \in \mathbb N_q
\end{align}
\end{subequations}
with $Q^i = \blkdiag(Q^i_1, \dots, Q^i_{m_i})$, $W^i = \begin{bmatrix} (W^i_1)^T & \dots & (W^i_{m_i})^T\end{bmatrix}^T$, $b^i = (b^i_1, \dots, b^i_{m_i})$, $x^i = \left(\svec(X^i_1), \dots, \svec(X^i_{m_i})\right)$ and $X^i = \blkdiag(X^i_1, \dots, X^i_{m_i})$. It is possible to solve this problem using a primal-dual interior-point method, \cite{wri:97}, \cite{tod:96}. Next we briefly discuss the main stages of such a method. The Karush-Kuhn-Tucker, KKT, optimality conditions for this problem are given as
\begin{subequations}\label{eq:KKT}
\begin{align}
\sum_{i = 1}^q \left((W^i)^Tv^i + (A^i)^T\bar v^i + (D^i)^T\lambda^i \right) &= -\sum_{i = 1}^q c^i, \\
(Q^i)^T v^i - z^i &= 0, \quad i \in \mathbb N_q,\\
X^iZ^i &= 0, \quad i \in \mathbb N_q, \label{eq:KKT-c}\\
\diag(\lambda^i)\left( D^iy-g^i \right) &= 0, \quad i \in \mathbb N_q,\label{eq:KKT-d}\\
Q^i x^i + W^i y &= b^i, \quad i \in \mathbb N_q, \\
A^i y &= \bar b^i, \quad i \in \mathbb N_q,
\end{align}
\end{subequations}
together with $D^i x \leq g^i$ and $X^i \succeq 0$, where $Z^i = \blkdiag(Z^i_1, \dots, Z^i_{m_i})$ and $z^i = \left(\svec(Z^i_1), \dots, \svec(Z^i_{m_i})\right)$. Any solution to this set of nonlinear equations is optimal for~\eqref{eq:GeneralSDP}. Within a primal-dual interior-point method, we set out to compute a solution to \eqref{eq:GeneralSDP}, by considering a sequence of perturbed KKT conditions where \eqref{eq:KKT-c} and \eqref{eq:KKT-d} are modified~as
\begin{align*}
X^iZ^i &= \delta I, \quad i \in \mathbb N_q,\\
\diag(\lambda^i)\left( D^iy-g^i \right) &= -\delta \mathbf 1, \quad i \in \mathbb N_q.
\end{align*}
where $\delta >0$ is the perturbation parameter. Particularly at each iteration, given feasible iterates $\lambda^i > 0$, $y$ so that $D^i y > g^i$ and $X^i \succ 0$ for $i = 1, \dots, q$, the primal-dual search directions are computed by solving a linearized version of the perturbed KKT conditions, given as
\small
\begin{subequations}\label{eq:KKTLin}
\begin{align}
\sum_{i=1}^q \left((W^i)^T \Delta v^i + (A^i)^T \Delta \bar v^i + (D^i)^T \Delta \lambda^i\right) &= r_{\textrm{d,lin}},\\
(Q^i)^T \Delta v^i - \Delta z^i &= r_{\textrm{d}}^i, \ i \in \mathbb N_q,\\
U^i \Delta x^i + F^i \Delta z^i &= r_{\textrm{c}}^i, \ i \in \mathbb N_q, \\
\diag(\Delta \lambda^i) (D^i y - g^i) + \diag(\lambda^i) D^i \Delta y &= r_{\textrm{c,lin}}^i, \ i \in \mathbb N_q,\\
Q^i \Delta x^i + W^i \Delta y &= r_{\textrm{p}}^i, \quad i \in \mathbb N_q, \\
A^i \Delta y &= r_{\textrm{p,lin}}^i, \quad i \in \mathbb N_q,
\end{align}
\end{subequations}
\normalsize
with $U^i = \blkdiag(U^i_1, \dots, U^i_{m_i})$, $F^i = \blkdiag(F^i_1, \dots, F^i_{m_i})$, where given
\begin{equation}\label{eq:scaling}
\begin{split}
W^i_j :&= (X^i_j)^{\frac{1}{2}} \left( (X^i_j)^{\frac{1}{2}} Z^i_j (X^i_j)^{\frac{1}{2}} \right)^{-\frac{1}{2}} (X^i_j)^{\frac{1}{2}}\\
& = (Z^i_j)^{-\frac{1}{2}} \left( (Z^i_j)^{\frac{1}{2}} X^i_j (Z^i_j)^{\frac{1}{2}} \right)^{\frac{1}{2}} (Z^i_j)^{-\frac{1}{2}},
\end{split}
\end{equation}
$W^i_j =: G^i_j(G^i_j)^T$ and $D^i_j = (G^i_j)^{-1}$, we have $U^i_j = D^i_j \otimes_s (D^i_j)^{-T}Z^i_j$ and $F^i_j = D^i_jX^i_j \otimes_s (D^i_j)^{-T}$. Furthermore, the residuals are given as
\begin{subequations}
\begin{align}
r_{\textrm{d,lin}} &= \sum_{i = 1}^q \underbrace{-c^i - (Q^i)^Tv^i - (A^i)^T\bar v^i - (D^i)^T\lambda^i}_{r_{d,lin}^i}\\
r_{\textrm{d}}^i &= z^i - (Q^i)^Tv^i, \quad i \in \mathbb N_q,\\
r_{\textrm{c}}^i &= \svec(\delta I - H_{D^i_j}(X^i_jZ^i_j)), \quad i \in \mathbb N_q,\\
r_{\textrm{c,lin}}^i &= -\delta \mathbf 1 - \diag(\lambda^i)(D^iy-g^i), \quad i \in \mathbb N_q,\\
r_{\textrm{p}}^i &= b^i - W^ix^i -Q^i y  , \quad i \in \mathbb N_q,\\
r_{\textrm{p,lin}}^i &= \bar b^i - A^i y  , \quad i \in \mathbb N_q,
\end{align}
\end{subequations}
where $H_D(M) = 1/2(DMD^{-1} + D^{-T}MD^{T})$. Having computed the search directions, suitable primal and dual step sizes, i.e., $t_d$ and $t_p$, are calculated so as to guarantee feasibility of the iterates with respect to inequality constraints and persistent reduction of residual norms, see e.g., \cite{tod:96} and references therein, which then allows us to update the iterates. This process is then repeated until certain stopping criteria are satisfied, which commonly depend on the residual norms and the size of the perturbation parameter. A generic description of a primal-dual interior-point method is given in Algorithm \ref{alg:PDIPM}.
 \begin{algorithm}[tb]
\caption{Primal-dual Interior-point Method}\label{alg:PDIPM}
\begin{algorithmic}[1]
\small
\State{Given feasible iterates with respect to inequality constraints}
\Repeat
\State{Compute the primal-dual search directions}
\State{Compute primal and dual step sizes}
\State{Update primal and dual iterates}
\State{Update the perturbation parameter}
\Until{stopping criteria is satisfied}
\normalsize
\end{algorithmic}
\end{algorithm}
The most computationally demanding step at every iteration of a primal-dual interior-point method, concerns the computation of the search directions. This requires solving the linear system of equations in \eqref{eq:KKTLin}, which can be written more compactly as
\begin{multline}\label{eq:KKTCompact}
\begin{bmatrix} W^T & A^T & & & & D^T \\ Q^T & & & & -I & \\ & & U & & F & \\ & & & \Lambda D & & E \\ & & Q & W & & \\& & & A & & \end{bmatrix} \begin{bmatrix}  \Delta v \\ \Delta \bar v \\ \Delta x \\ \Delta y \\ \Delta z \\ \Delta \lambda \end{bmatrix} = \begin{bmatrix} r_{\textrm{d,lin}}\\ r_{\textrm{d}} \\ r_{\textrm{c}} \\ r_{\textrm{c,lin}}\\ r_{\textrm{p}} \\ r_{\textrm{p,lin}}  \end{bmatrix}
\end{multline}
where
\begin{align*}
W &= \begin{bmatrix} (W^1)^T & \dots & (W^q)^T\end{bmatrix}^T,\\
A &= \begin{bmatrix} (A^1)^T & \dots & (A^q)^T\end{bmatrix}^T,\\
Q &= \blkdiag(Q^1, \dots, Q^q),\\
U &= \blkdiag(U^1, \dots, U^q),\\
F &= \blkdiag(F^1, \dots, F^q),\\
D &= \begin{bmatrix} (D^1)^T & \dots & (D^q)^T\end{bmatrix}^T,\\
\Lambda &= \blkdiag(\Lambda^1, \dots, \Lambda^q), \quad \Lambda^i = \diag(\lambda^i),\\
E &= \blkdiag(E^1, \dots, E^q), \quad E^i = \diag(D^iy - g^i),
\end{align*}
and the variables and the right hand side terms correspond to all variables and residuals stacked. One way to solve this system of equations is by first eliminating the third and fourth row equations as
\begin{subequations}
\begin{align}
\Delta z &= F^{-1}\left( r_c - U \Delta x \right), \\
\Delta \lambda &= E^{-1}\left( r_{\textrm{c,lin}} - \Lambda D \Delta y  \right),
\end{align}
\end{subequations}
which is possible since $F$ and $E$ are both invertible, see e.g., \cite{tod:96, wri:97}. This then allows us to rewrite \eqref{eq:KKTCompact} as
\begin{multline}
\begin{bmatrix} -D^TE^{-1}D &  & W^T & A^T \\  & F^{-1}U & Q^T & \\ W & Q &  & \\A & & &  \end{bmatrix} \begin{bmatrix}  \Delta y \\ \Delta x \\ \Delta v \\ \Delta \bar v \end{bmatrix} = \begin{bmatrix} r_{\textrm{lin}}\\ r \\  r_{\textrm{p}} \\ r_{\textrm{p,lin}}  \end{bmatrix}
\end{multline}
where $r_{\textrm{lin}} = r_{\textrm{d,lin}} - D^T E^{-1} r_{\textrm{c,lin}}$ and $r = r_{\textrm{d}} - F^{-1} r_{\textrm{c}}$. Notice that this set of linear equations also defines the optimality conditions for the convex quadratic program (QP)
\small
\begin{subequations}\label{eq:KKTQP}
\begin{align}
\minimize \quad& \begin{bmatrix}\Delta y \\ \Delta x \end{bmatrix}^T \begin{bmatrix} -D^TE^{-1}D & \\ & F^{-1}U \end{bmatrix} \begin{bmatrix} \Delta y \\ \Delta x  \end{bmatrix} -\notag\\ &\hspace{50mm} \begin{bmatrix} r_{\textrm{lin}}\\ r \end{bmatrix}^T \begin{bmatrix} \Delta y \\ \Delta x  \end{bmatrix}\\
\subject \quad& W\Delta y + Q \Delta x = r_{\textrm{p}} \\
&  A \Delta y = r_{\textrm{p,lin}}
\end{align}
\end{subequations}
\normalsize
For the localization problem, this QP has a particular structure which enables us to solve it distributedly and efficiently, using message-passing. Next we briefly discuss this algorithm for the sake of completeness and to provide a better understanding of the presented material.

\subsection{Solving Coupled Optimization Problems Using Message-passing}
\begin{figure}[t]
\begin{center}
\includegraphics[width=3cm]{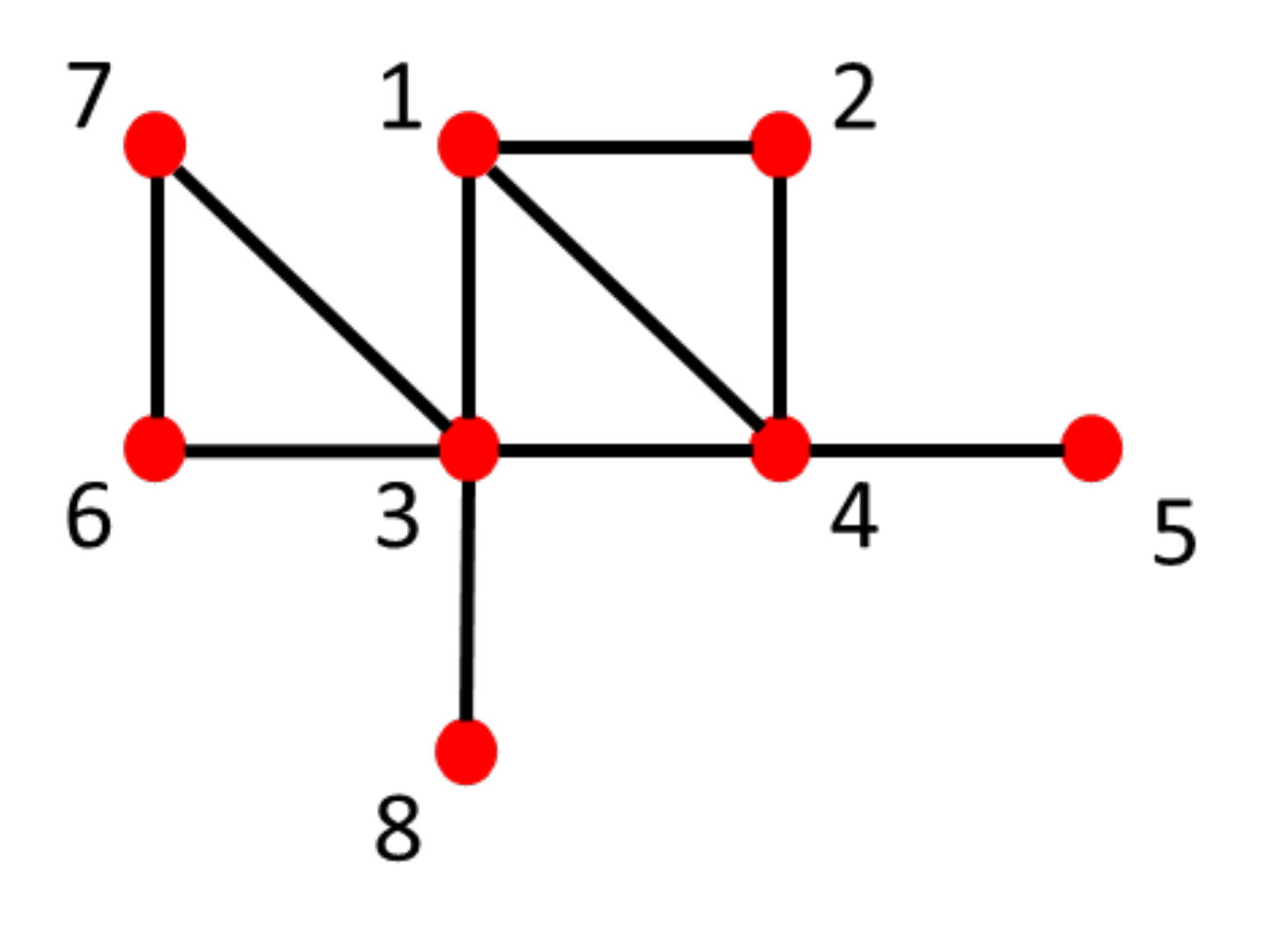}    
\caption{ The sparsity graph for the problem in~\eqref{eq:Example}.}
\label{fig:Example}
\vspace{-5mm}
\end{center}
\end{figure}

Consider the following coupled optimization problem
\begin{align}\label{eq:UncostrainedProblem}
\minimize & \quad F_1(x) + F_2(x) + \dots + F_q(x),
\end{align}
where $x \in \mathbb R^n$ and and the functions $F_i \ : \ \mathbb R^n \rightarrow \mathbb R$ for $i \in \mathbb N_q$ are convex. Also we assume that each term in the objective function (each subproblem) only depends on a few variables. Let us denote the indices of the variables that appear in the $i$th term, $F_i$, by $J_i$. This definition allows us to rewrite the problem in \eqref{eq:UncostrainedProblem} as
\begin{align}\label{eq:UncostrainedProblemReform}
\minimize & \quad \bar F_1(x_{_{J_1}}) +\bar F_2(x_{_{J_2}}) + \dots +\bar F_q(x_{_{J_q}}),
\end{align}
where $x_{_{J_i}} = E_{J_i} x$. The functions $\bar F_i \ : \ \mathbb R^{|J_i|} \rightarrow \mathbb R$ are lower dimensional descriptions of $F_i$s such that $F_i(x) = \bar F_i(E_{J_i}x)$ for all $x \in \mathbb R^n$ and $i \in \mathbb N_N$. We also define $\mathcal I_j$ as the set of indices of terms in the cost function that depend on $x_j$, i.e., $\{ i \ | \ j \in J_i \}$. The sets $J_i$ for $i \in \mathbb N_q$ and $\mathcal I_j$ for $j \in \mathbb N_n$ provide a clear mathematical description of the coupling structure in the problem. It is also possible to describe the coupling structure in the problem graphically, using graphs. For this purpose, we introduce the sparsity graph. The \emph{sparsity graph} $G_s(V_s, \mathcal E_s)$ of a coupled problem is an undirected graph with the vertex set $V_s = \left\{ 1, \dots, n\right\}$ and the edge set $\mathcal E_s = \left\{ (i, j) \ | \ i, j \in V_s, \ \mathcal I_i \cap \mathcal I_j \neq \emptyset \right\}$. As an example consider the following problem
\begin{multline}\label{eq:Example}
\minimize_x \quad \bar F_1(x_1, x_3, x_4) + \bar F_2(x_1, x_2, x_4) + \bar F_3(x_4, x_5) +\\ \bar F_4(x_3, x_6, x_7) + \bar F_5(x_3, x_8).
\end{multline}
The sparsity graph for this problem are illustrated in Figure~\ref{fig:Example}.
It is possible to devise scalable or distributed algorithms for solving the problem in \eqref{eq:UncostrainedProblem}. In this paper we focus on message-passing.

Consider the problem in \eqref{eq:UncostrainedProblemReform}, and assume that its sparsity graph is chordal. Let its set of cliques be given as $\mathbf C_{G_s} = \{ J_1, \dots, J_q\}$ and $T_s(V_t, \mathcal E_t)$ be a clique tree over the cliques. It is possible to solve the problem in \eqref{eq:UncostrainedProblemReform} distributedly, using an algorithm with the clique tree as its computational graph. That means each node in the tree corresponds to a computational agent and they communicate/collaborate with one another if there is an edge between them. Recall that each node in the clique tree is assigned a clique of the sparsity graph, i.e., $J_i$. In such a setting, we also assign each term in the objective function (each subproblem), i.e., $\bar F_i$, to each agent $i$. We can now describe how the problem in \eqref{eq:UncostrainedProblemReform} can be solved using message-passing by performing an upward-downward pass through the clique tree. The message-passing algorithm starts from the agents at the leaves of the tree, i.e., all $i \in \leaves(T)$, where every such agent computes the following message
\begin{align}\label{eq:mijLeaves}
m_{i\parent(i)}(x_{_{S_{i\parent(i)}}}) = \minimum_{x_{_{R_{i\parent(i)}}}} \left\{  \bar F_i(x_{_{J_i}})\right\},
\end{align}
with $S_{i\parent(i)} := J_i \cap J_{\parent(i)}$ and $R_{i\parent(i)} := J_i \setminus S_{i\parent(i)}$ are the so-called separators and residuals, respectively, and communicates it to its corresponding parent, denoted by $\parent(i)$. Notice that this message is a functional and not a scalar value, and hence agent $i$ needs to communicate the functional form. Then every parent $j$ that has received these messages from its children, denoted by $\children(j)$, computes its corresponding message to its parent as

\begin{figure}[t]
\begin{center}
\includegraphics[width=9cm]{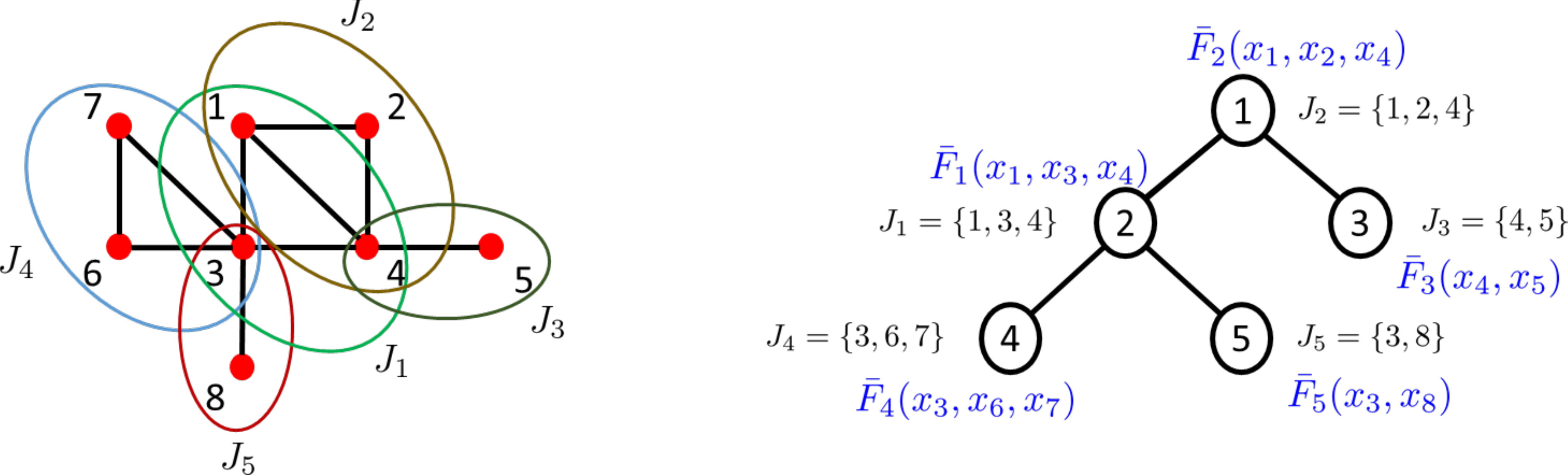}    
\caption{ The coupling and sparsity graphs for the problem in~\eqref{eq:Example}, illustrated on the right and left figures, respectively.}
\label{fig:ExampleClique}
\vspace{-5mm}
\end{center}
\end{figure}

\small
\begin{align}\label{eq:mij}
m_{j\parent(j)}(x_{_{S_{j\parent(j)}}}) = \minimum_{x_{_{R_{j\parent(j)}}}} \left\{   \bar F_j(x_{_{J_j}}) +  \sum_{k \in \children(j)} m_{kj}(x_{_{S_{kj}}}) \right\}.
\end{align}
\normalsize
This procedure is then continued until we arrive at the agent at the root. At this point, the agent at the root, indexed $r$, having received all messages from its children can compute the optimal solution for its corresponding variables specified by $J_r$ as
\begin{align}\label{eq:RLocalProblem}
x^\ast_{_{J_r}} = \argmin_{x_{_{J_r}}} \left\{  \bar F_k(x_{_{J_r}})  + \sum_{k \in \children(r)} m_{kr}(x_{_{S_{rk}}}) \right\}.
\end{align}
This agent then having computed its optimal solution, communicates this solution to its children, at which point every such agent $i \in \children(r)$ computes its optimal solution as
\begin{multline}\label{eq:LocalProblempar}
x^\ast_{_{J_i}} = \argmin_{x_{_{J_i}}} \left\{  \bar F_i(x_{_{J_i}})   + \sum_{k \in \children(i)} m_{ki}(x_{_{S_{ik}}}) + \right. \\ \left. \frac{1}{2} \left\| x_{_{S_{\parent(i)i}}} - \left( x_{_{S_{\parent(i)i}}}^\ast \right)^{\parent(i)} \right\|^2 \right\},
\end{multline}
where $\left( x_{_{S_{\parent(i)i}}}^\ast \right)^{\parent(i)}$ is the the computed optimal solution by the parent $\parent(i)$. This procedure is continued until we reach the agents at the leaves. At this point all agents have computed their corresponding optimal solution and the algorithm can be terminated, and hence, we have convergence after one upward-downward pass through the tree, \cite{kho:15c}, \cite{kol:09}. Let us now illustrate this procedure using an example. Consider the example given in \eqref{eq:Example}. The sparsity graph of this problem is chordal and its cliques are marked in Figure \ref{fig:ExampleClique} on the left. A clique tree for this graph is illustrated in the same figure on the right, where also a valid subproblem assignment is presented.
As was discussed above we start the message-passing from the leaves of the tree, particularly agents 3, 4 and 5. These agents compute and communicate their messages to their corresponding parents as
\begin{align*}
m_{32}(x_4) &= \minimum_{x_5}\left\{ \bar F_3(x_4, x_5)\right\}\\
m_{41}(x_3) & =   \minimum_{x_6, x_7}\left\{\bar F_4(x_3, x_6, x_7)\right\} \\
m_{51}(x_3) & =  \minimum_{x_8}\left\{ \bar F_5(x_3, x_8)\right\}.
\end{align*}
At this point agent 2 has received all messages from its children and can in turn compute and communicate its message to its parent as
\begin{align*}
m_{12}(x_1, x_4) &= \minimum_{x_3}\left\{ m_{41}(x_3) + m_{51}(x_3) + \bar F_1(x_1, x_3, x_4)\right\}.
\end{align*}
This completes the upward pass and now the agent at the root, i.e., agent 2, can compute its optimal solution as
\begin{multline*}
(x_1^*, x_2^*, x_4^*) = \argmin_{x_1, x_2, x_4} \\\quad\left\{ m_{12}(x_1, x_4) + m_{32}(x_4) + \bar F_2(x_1, x_2, x_4)\right\},
\end{multline*}
which initiates the downward pass. Agent 2 will then communicate $x_1^*, x_4^*$ and $x_4^*$ to agents 2 and 3 respectively, where they compute their corresponding optimal solution for the remainder of their variables as
\begin{align*}
x_3^* &= \argmin_{x_3} \quad\left\{ m_{41}(x_3) + m_{51}(x_3) + \bar F_1(x_1^*, x_3, x_4^*)\right\}\\
x_5^* &= \argmin_{x_5} \quad\left\{ \bar F_3(x_4^*, x_5)\right\}.
\end{align*}
The last step of the downward pass is then accomplished by agent 2 communicating $x_3^*$ to agents 4 and 5, and these agents computing their optimal solution as
\begin{align*}
(x_6^*, x_7^*) &= \argmin_{x_6, x_7} \quad\left\{ \bar F_4( x_3^*, x_6, x_7)\right\}\\
x_8^* &= \argmin_{x_8} \quad\left\{ \bar F_5(x_3^*, x_8)\right\},
\end{align*}
which finishes the algorithm. Notice that the message-passing algorithm described in this section can be viewed as dynamic programming over trees. Next we discuss how message-passing can be used within the primal-dual method.

\subsection{Distributed Computations In Primal-dual methods}
The problem in \eqref{eq:KKTQP} can be written as

\small
\begin{subequations}\label{eq:KKTQPSum}
\begin{align}
\minimize \quad& \sum_{i = 1}^q \begin{bmatrix}\Delta y \\ \Delta x^i \end{bmatrix}^T \begin{bmatrix} H^i & \\ & (F^i)^{-1}U^i \end{bmatrix} \times \notag\\ & \hspace{35mm} \begin{bmatrix} \Delta y \\ \Delta x^i  \end{bmatrix} - \begin{bmatrix} r^i_{\textrm{lin}}\\ r^i \end{bmatrix}^T \begin{bmatrix} \Delta y \\ \Delta x^i  \end{bmatrix}\\
\subject &\notag \\\quad& \begin{rcases*} W^i\Delta y + Q^i \Delta x^i = r^i_{\textrm{p}}, \\
  A^i \Delta y = r^i_{\textrm{p,lin}}, \end{rcases*}\quad i \in \mathbf N_q
\end{align}
\end{subequations}
\normalsize
where $H^i = -(D^i)^T(E^i)^{-1}D^i$ and $r^i_{\textrm{lin}} = r^i_{\textrm{d,lin}} - (D^i)^T (E^i)^{-1} r^i_{\textrm{c,lin}}$ and $r^i = r^i_{\textrm{d}} - (F^i)^{-1} r^i_{\textrm{c}}$. This problem can be viewed as a combination of $q$ subproblems, where each of which is defined by a term in the objective function and its corresponding equality constraints. Notice that the coupling among the subproblems does not stem from the matrix variables and on the surface all subproblems seem to be coupled to one another through the linear variables directions $\Delta y$. However, for the localization problem in \eqref{eq:MLConstrainedScattered3}, each subproblem only relies on a certain elements of $\Delta y$. This can be seen by first noticing that the linear variables for each subproblem $k$ is given by $\vectri(S_{C_kC_k})$, $\Lambda_{ij}, D_{ij}$  for $(i,j) \in \phi_k$ and $x^i_s, \Xi_{ij} , Z_{ij}$ for $j \in \textrm{Ne}_a(i)$ and $i \in \bar \phi_k$. Let us assume that the indices of elements of $\Delta y$ that correspond to these variables be given by set $J_k$. We can then rewrite the problem in \eqref{eq:KKTQPSum} as

\small
\begin{subequations}\label{eq:KKTQPSum}
\begin{align}
\minimize \quad& \sum_{i = 1}^q \begin{bmatrix}\Delta y_{_{J_i}} \\ \Delta x^i \end{bmatrix}^T \begin{bmatrix} \bar H^i & \\ & (F^i)^{-1}U^i \end{bmatrix} \times \notag\\ & \hspace{35mm} \begin{bmatrix} \Delta y_{_{J_i}} \\ \Delta x^i  \end{bmatrix} - \begin{bmatrix}\bar r^i_{\textrm{lin}}\\ r^i \end{bmatrix}^T \begin{bmatrix} \Delta y_{_{J_i}} \\ \Delta x^i  \end{bmatrix}\\
\subject \quad& \bar W^i\Delta y_{_{J_i}} + Q^i \Delta x^i = r^i_{\textrm{p}}, \quad i = 1, \dots, q \\
&  \bar A^i \Delta y_{_{J_i}} = r^i_{\textrm{p,lin}}, \quad i = 1, \dots, q
\end{align}
\end{subequations}
\normalsize
where $\bar H^i = E_{J_i}H^iE_{J_i}^T$, $\bar r^i_{\textrm{lin}} = E_{J_i}r^i_{\textrm{lin}}$, $\bar A^i = A^iE_{J_i}^T$ and $\bar W^i = W^i E_{J_i}^T$. Through the use of indicator functions, this problem can be written as
\small
\begin{align}
\minimize \quad& \sum_{i = 1}^q \begin{bmatrix}\Delta y_{_{J_i}} \\ \Delta x^i \end{bmatrix}^T \begin{bmatrix} \bar H^i & \\ & (F^i)^{-1}U^i \end{bmatrix}   \begin{bmatrix} \Delta y_{_{J_i}} \\ \Delta x^i  \end{bmatrix} - \notag\\ & \hspace{25mm}\begin{bmatrix}\bar r^i_{\textrm{lin}}\\ r^i \end{bmatrix}^T \begin{bmatrix} \Delta y_{_{J_i}} \\ \Delta x^i  \end{bmatrix} + \mathcal I_{\mathcal C_i}(\Delta y_{_{J_i}}, \Delta x^i)
\end{align}
\normalsize
where $\mathcal C_i = \{ (\Delta y_{_{J_i}},\Delta x^i)\ | \  \bar W^i\Delta y_{_{J_i}} + Q^i \Delta x^i = r^i_{\textrm{p}}, \ \bar A^i \Delta y_{_{J_i}} = r^i_{\textrm{p,lin}}\}$ and
\begin{align*}
\mathcal I_{\mathcal C_i}\left(x\right) = \begin{cases} 0 \hspace{6mm} x \in \mathcal C_i  \\ \infty \hspace{4mm} \text{Otherwise} \end{cases}
\end{align*}
This problem is in the same format as \eqref{eq:UncostrainedProblemReform}. It is now possible to see that the coupling comes from the fact that for some $C_i$ and $C_j$, $C_i \cap C_j \neq \emptyset$. Recall that one way to describe the intersection among the cliques of the inter-sensor measurement graph can be described using its clique tree, $T(V_t, \mathcal E_t)$. The sparsity graph of this problem is in fact chordal with cliques defined by the variables that appear in each subproblem. Furthermore, the clique tree for the sparsity graph of this problem has the same structure as that of the inter-sensor measurement graph. This is the case since the ordering defined by this tree defines perfect elimination ordering for the sparsity graph, see \cite{gol:04} for more details. Consequently, this problem can be solved distributedly using message-passing as discussed above. As a result, we can compute the primal-dual search directions for the problem in \eqref{eq:MLConstrainedScattered3} distributedly, by an upward-downward pass through the clique tree. Notice that the messages for solving this problem are quadratic functions, and hence the hessian and linear term that describes this function need to be communicated. The remaining stages of a primal-dual interior-point method can also be done distributedly over the clique tree. For the sake of brevity, we here do not discuss the details any further, for more info see, \cite{kho:15c} and \cite{kho:15d}. A summary of our proposed distributed localization method is given in Algorithm \ref{alg:Local}.
\begin{algorithm}[tb]
\caption{Distributed Primal-dual Localization Algorithm, DPDLA}\label{alg:Local}
\begin{algorithmic}[1]
\small
\State{Given the inter-sensor measurement graph $G_r(V_r, \mathcal E_r)$, its cliques set $C_{G_r} = \{ C_1, \dots, C_q \}$ and a clique tree over its cliques $T(V_t, \mathcal E_t)$ with $V_t = \{ 1, \dots, q\}$}
\State{Conduct assignments such that the assignment rules in Section \ref{sec:decompositionloc} are satisfied, for instance using Algorithm \ref{alg:Ass}}
\State{Each agent $i \in \mathbb N_q$ forms its corresponding subproblem}
\State{Given feasible initial primal and dual iterates with respect to inequality constraints}
\Repeat
\State{Compute the primal-dual search directions distributedly using message-passing over $T(V_t, \mathcal E_t)$}
\State{Compute primal and dual step sizes distributedly (this can be done by performing an upward-downward pass through $T(V_t, \mathcal E_t)$, see \cite[Sec. 6.4]{kho:15c},\cite[Sec. V-B]{kho:15d})}
\State{Update primal and dual iterates}
\State{Update the perturbation parameter and the compute the stopping criteria distributedly (this can be done by performing an upward-downward pass through $T(V_t, \mathcal E_t)$, see \cite[Sec. 6.4]{kho:15c},\cite[Sec. V-B]{kho:15d})}
\Until{stopping criteria is satisfied}
\normalsize
\end{algorithmic}
\end{algorithm}

\subsection{Computational and Communication Complexity}

At each iteration of the primal-dual method, we need to conduct three upward-downward passes, namely one for computing the primal-dual directions, one for computing the primal and dual step sizes and one for updating the perturbation parameter and checking the termination condition. This means that if the primal-dual method converges within $p$ iterations, the algorithm converges within $3\times2\times p\times h$ steps where $h$ is the height of the considered clique tree. Furthermore, during the execution of the algorithm, each agent is required to communicate twice with its neighbors during each upward-downward pass. Once with its parent during the upward pass and once with its children during the downward pass. Consequently, the total number of times each agent needs to communicate with its neighbors is given by $3\times 2\times p$.

Among the upward and downward passes, the upward pass for computing the search directions, is the most computationally demanding and communication intensive one. Particularly, during this upward pass each agent needs to compute a factorization of a relatively small matrix to compute its message to the parent, see \cite[Sec. 6.2]{kho:15c}. This needs to be done once at every primal-dual iteration, which means that in total each agent is required to compute $p$ factorizations during the run of Algorithm \ref{alg:Local}. Also recall that during these upward passes, each agent needs to communicate a quadratic functional to its parent. This entails sending the data matrices that define the quadratic function. Depending on the number of variables shared between each agent and its parent, the information that needs to be communicated can be considerable. Notice that the computational burden of the other upward and downward passes are comparatively trivial. Moreover, the information that needs to be communicated during these upward and downward passes is limited to a few scalars. Due to this, in the remainder of this section, we discuss the computational and communication burden for each agent during the upward pass for computing the search directions.

Firstly, recall that each subproblem $k$ in \eqref{eq:KKTQPSum}, depends on variables
\begin{align*}
&\vectri(S_{C_kC_k}),\\ &\Lambda_{ij}, D_{ij} \quad \textrm{for} \ (i,j) \in \phi_k, \\ &x^i_s, \Xi_{ij} , Z_{ij} \quad \textrm{for} \  j \in \textrm{Ne}_a(i), \ i \in \bar \phi_k, \\ & T^k, \\ & \Gamma^{ij} \quad \textrm{for} \ (i,j) \in \phi_k, \\ &\Phi^{ij} \quad \textrm{for} \  j \in \textrm{Ne}_a(i), \ i \in \bar \phi_k.
\end{align*}
Let us assume that each agent $k$ is assigned $b_k$ and $a_k$ inter-sensor and anchor-sensor range measurements, respectively. The number of variables that appear in each subproblem $k$ is then given as

\small
\begin{multline}
n_k = \frac{|C_k|(|C_k|+1)}{2} + 2b_k + 2|C_k| + 2a_k + \\ \frac{(|C_k|+2)(|C_k|+3)}{2} + 3b_k + 3a_k.
\end{multline}
\normalsize
Notice that the number of equality constraints defined by each range measurement is equal to four, see \eqref{eq:MLConstrainedScattered2-b} and \eqref{eq:MLConstrainedScattered2-c}. Consequently the number of equality constraints for each subproblem $k$ is given as

\small
\begin{align}
e_k = 4b_k + 4a_k + \frac{(|C_k|+2)(|C_k|+3)}{2}.
\end{align}
\normalsize
Let us define $U_k = C_k \cap C_{\parent(k)}$. The variables that are shared between agent $k$ and its parent are given as $x^i_s$ for $i \in U_k$ and $\vectri(S_{U_kU_k})$. The number of these variables is then $s_k = 2|U_k| + |U_k|(|U_k|+1)/2$. The number of variables that agent $k$ does not share with its parent is then $r_k = n_k - s_k$. Each agent in order to compute the message to its parent, needs to factorize a symmetric indefinite matrix, see \cite[Sec. 6.2]{kho:15c}. The size of this matrix depends on the number of equality constraints for its subproblem and the variables it does not share with its parent. Hence, the size of this matrix is given by $r_k + e_k$. Moreover recall that the messages are quadratic functions of the variables that are shared between two agents. Consequently, each agent in order to communicate this functional to its parent would need to send $s_k(s_k+1)/2+s_k$ scalars to its parent. We can now summarize the dominant computational and communication burden for each agent with the following items.
\begin{itemize}
\item The size of the matrix that needs to be factorized by each agent $k$ grows quadratically with the number of sensors assigned to the agent and linearly with the number of range measurements assigned to it. This number is also reduced quadratically with the number of variables that this agent shares with its parent.
\item The size of the information that each agent needs to communicate to its parent grows quadratically with the number of variables it shares with the parent.
\end{itemize}
\begin{rem}
Notice that these summarizing items also provide guidelines on how to devise heuristics to perform a better clustering of sensors. They also enable us to propose improvements to the measurement assignment strategy, in order to distribute the computations among agents in a more balanced manner. Despite this, for the sake of brevity and simplicity, such heuristics are not considered in this study.
\end{rem}
Next we investigate the performance of our proposed algorithm, using two sets of numerical experiments.
%
%
%

\section{Numerical Experiments}\label{sec:numerical}
In this section we compare the performance of our proposed distributed algorithm with that of presented in \cite{soa:15}. We refer to this algorithm as distributed disk relaxation algorithm (DDRA). To this end, we conduct two sets of experiments, one that relies on simulated data and one that is based on real data from \cite{pat:03}. Notice that we do not conduct a comparison with other algorithms, since a thorough comparison with DDRA has been conducted in \cite{soa:15}, which illustrated the superiority of their proposed algorithm to high performance algorithms in \cite{gho:13} and \cite{sim:14} both in accuracy and number of communications among agents.

\subsection{Experiments Using Simulated Data}\label{sec:numerical-A}

Our experiments based on simulated data concern networks of sensors with connected inter-sensor measurement graphs. In all experiments there are 9 anchors in the network which are uniformly distributed in the area. The experiments in this section are divided into two setups. In both setups, we consider a network of several sensors which are placed in a two-dimensional area, with their locations randomly generated using a uniform distribution. The noisy range measurements are generated as
\begin{align*}
\mathcal R_{ij} &= \left|\| (x^*_s)^i - (x^*_s)^j \|_2 + E_{ij}\right|, \quad  j \in \text{Ne}_r(i),\\
\mathcal Y_{ij} &= \left|\| (x^*_s)^i - (x_a)^j \|_2 + V_{ij}\right|, \quad  j \in \text{Ne}_a(i),
\end{align*}
where $(x^*_s)^i$ denotes the true location of the $i$th sensor. Furthermore we assume that all noises are gaussian and mutually independent, see also \cite{soa:15}. In the first setup we conduct experiments using a network 50 sensors in a $0.8 \times 0.8$ area. We consider four different measurement noise standard deviations, namely $0.01, 0.05$, $0.1$ and $0.3$, and for each noise level we generate 50 problem instances. In order to ensure that the generated inter-sensor measurement graph is loosely connected, we assume there exist a measurement between two sensors or between a sensor and an anchor if the distance between them is less than the communication range $r_c =0.2$.
\begin{figure}[t]
\begin{center}
\includegraphics[width=6.5cm]{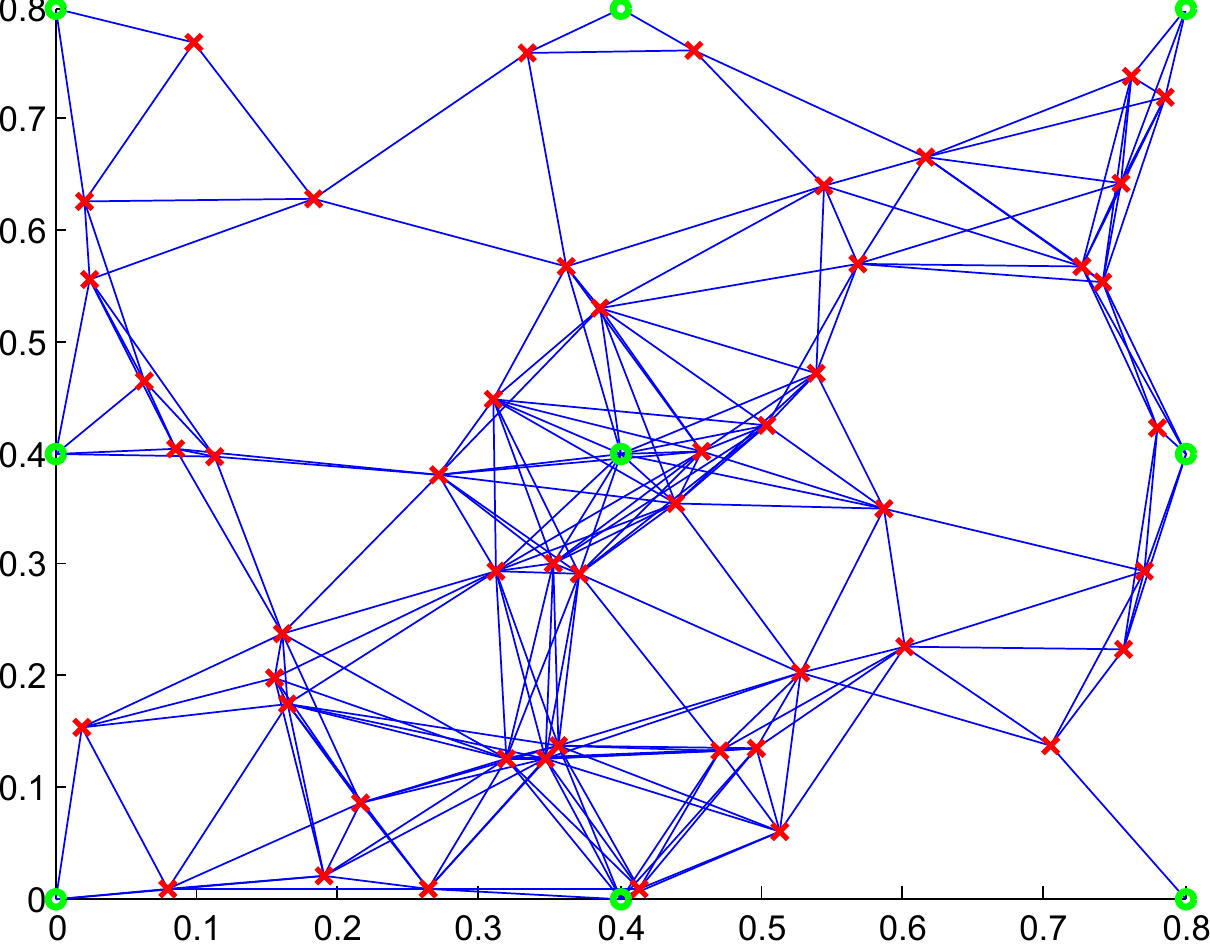}    
\caption{\footnotesize The sensor network considered for our experiment. Each red cross depicts one of the 50 sensors in the network and each green circle marks one of the 9 anchors. An edge between two nodes, implies existence of a range measurement between the two nodes.\normalsize}
\label{fig:network}
\vspace{-5mm}
\end{center}
\end{figure}
The resulting sensor network is depicted in Figure \ref{fig:network}. In this figure, the sensor nodes are marked with red crosses and the anchors are marked with green circles. As can be seen from the figure the inter-sensor measurement graph is connected. The performance of distributed algorithms are quantified using three measures. Namely (i) their accuracy based on the root mean squared error (RMSE) defined as
\begin{align}
\textrm{RMSE} = \sqrt{\frac{1}{MN}\sum_{i=1}^M\sum_{j=1}^N \| (x^*_s)^j - x_s^j(m) \|^2}
\end{align}
where $M$ is the number of experiments and the argument $m$ marks the computed estimate for the $m$th experiment, (ii) number of required iterations and communications to converge to a solution with a given accuracy and (iii) the computational time. Notice that both algorithms are run in a centralized manner. The algorithm in \cite{soa:15} is terminated if the norm of the gradient of its considered cost function is below $10^{-6}$. This threshold was chosen based on the authors experience, so as to guarantee DDRA generates accurate enough solutions. Figures~\ref{fig:ex11}--\ref{fig:ex13} illustrate the achieved results. In these figures and the ones to come the $*$-marked curves illustrate the results from DPDLA, whereas the o-marked curves show the results from DDRA.
\begin{figure}[t]
\begin{center}
\includegraphics[width=7.5cm]{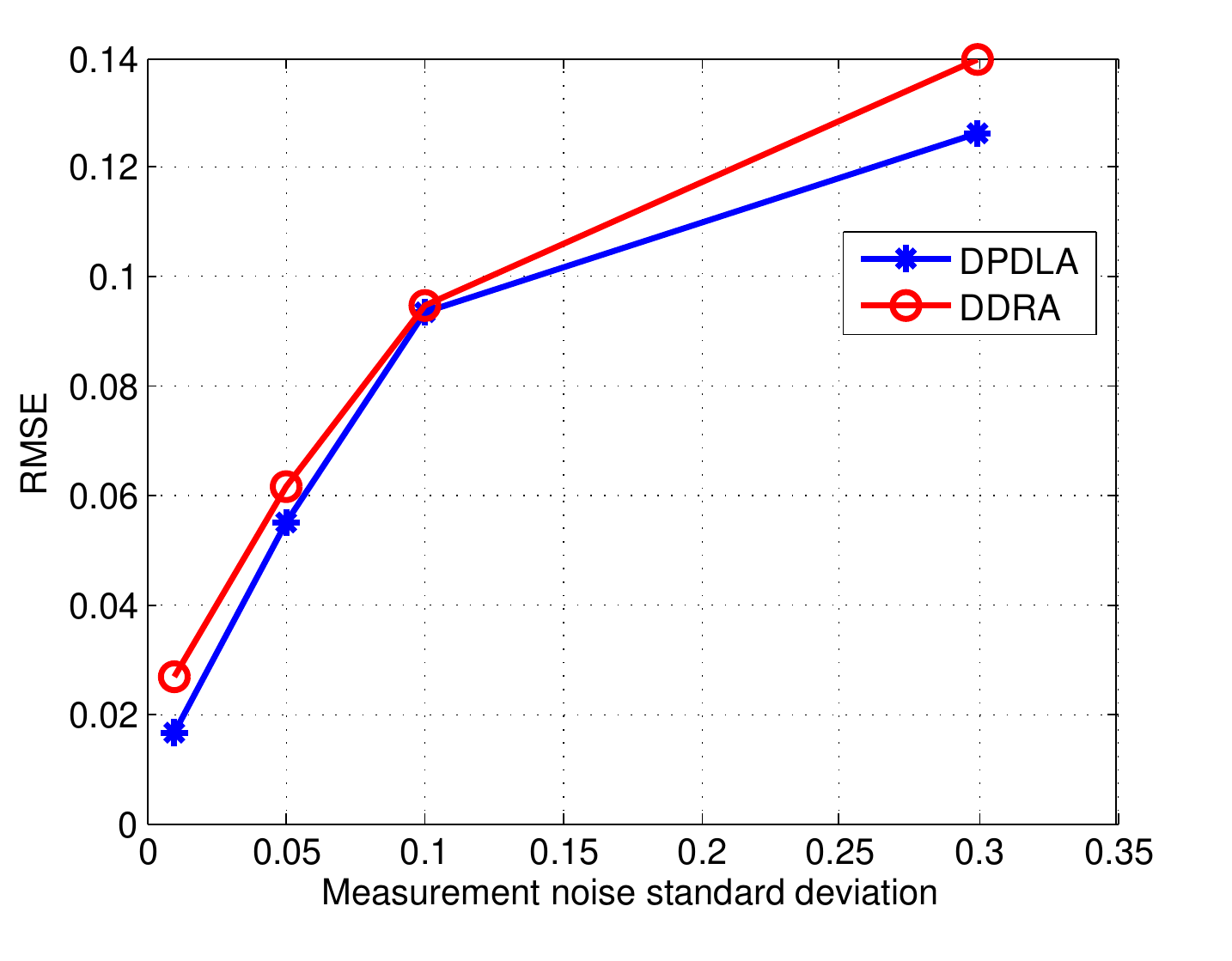}    
\caption{\footnotesize The RMSE results from the considered algorithms when applied to a network of 50 sensors, depicted in Figure \ref{fig:network}, with four different measurement noise standard deviation, namely $0.01, 0.05$, $0.1$ and $0.3$.\normalsize}
\label{fig:ex11}
\vspace{-5mm}
\end{center}
\end{figure}
As can be seen from Figure \ref{fig:ex11}, DPDLA outperforms or provides comparable accuracy with respect to DDRA for different levels of measurement noise. This shows the superiority of semidefinite relaxation to disk relaxation. The number of communications that each agent is required to conduct for each algorithm to converge to a solution is depicted in Figure~\ref{fig:ex12}. For these experiments, the considered clique tree for the inter-sensor measurement graph in Figure~\ref{fig:network}, has height~8, and the primal-dual method converged within around 10 iterations. As can be seen from the figure, DPDLA requires roughly two orders of magnitude less number of communications for computing a solution. The shaded areas depict the maximum and minimum values within the 50 instances for each of these quantities. Notice that this area for the results corresponding to DPDLA is not even visible. We can hence deduce that in comparison DDRA, the number of communications for DPDLA seems to be much less sensitive to the noise level and also to data realizations.
\begin{figure}[t]
\begin{center}
\includegraphics[width=7.5cm]{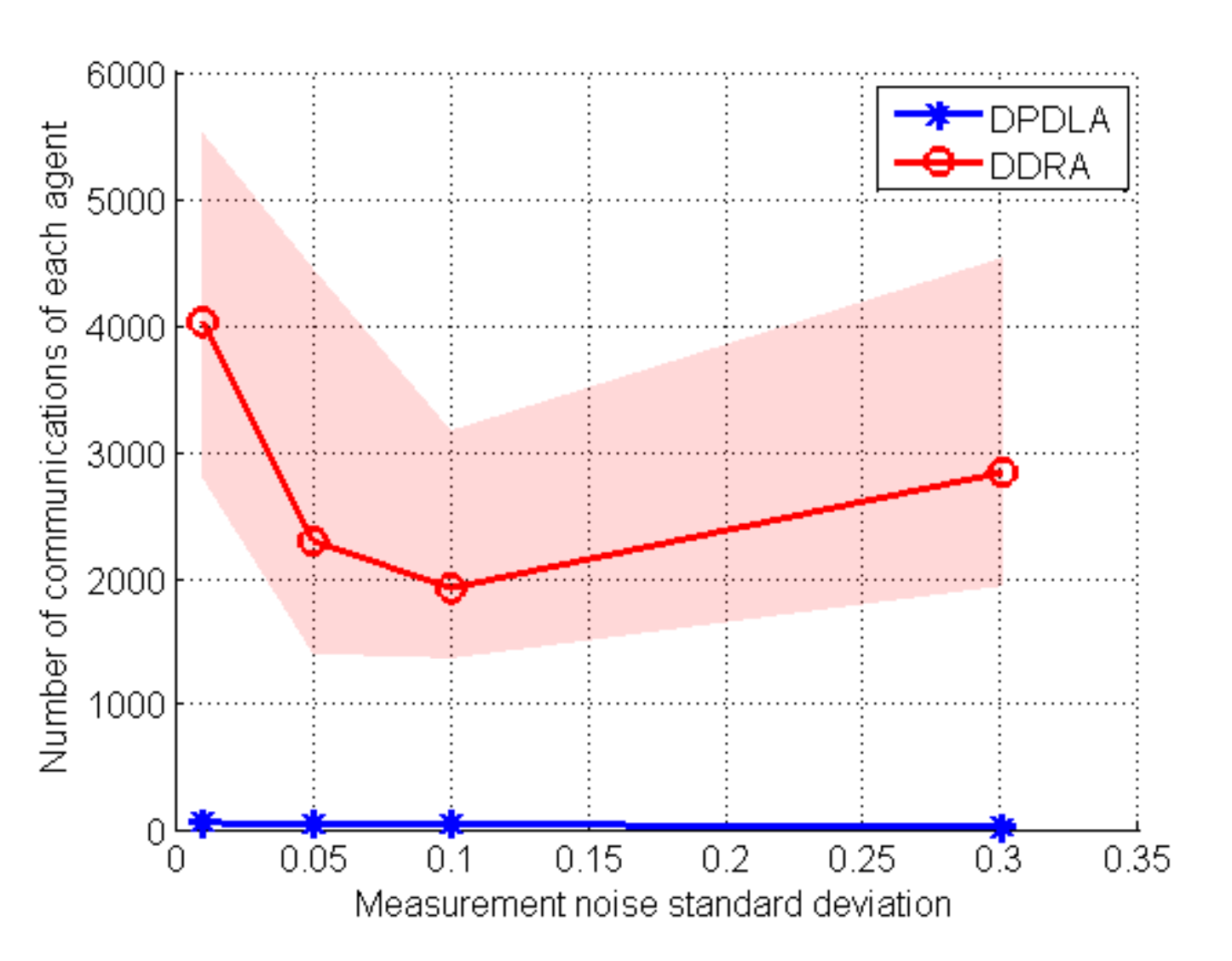}    
\caption{\footnotesize The number of communications that each agent needs to conduct for each of the algorithms to converge to a solution. The sensor network consists of 50 sensors, depicted in Figure \ref{fig:network}, with three different measurement noise standard deviation, namely $0.01, 0.05$, $0.1$ and~$0.3$. \normalsize}
\label{fig:ex12}
\vspace{-5mm}
\end{center}
\end{figure}
\begin{figure}[t]
\begin{center}
\includegraphics[width=7.5cm]{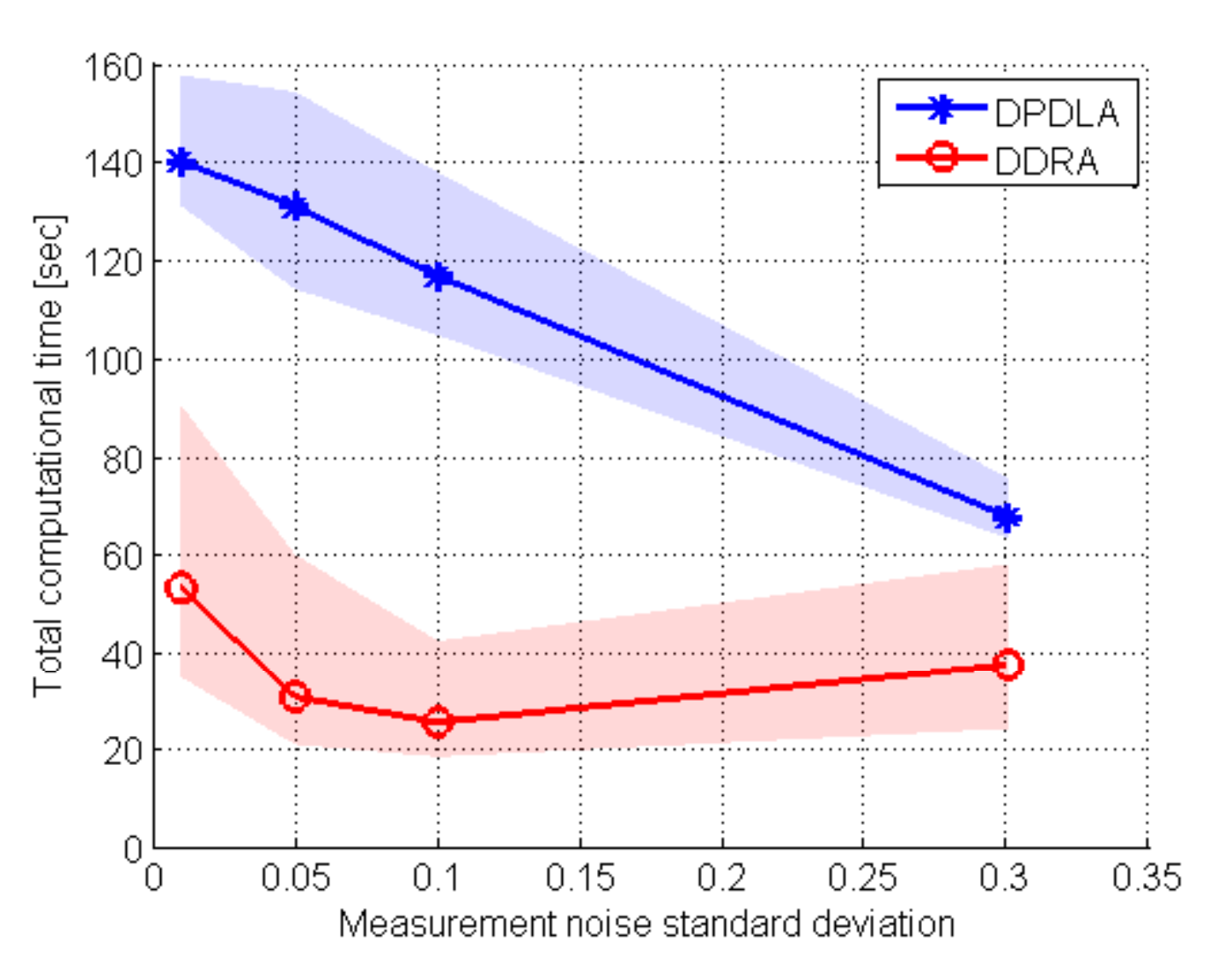}    
\caption{\footnotesize The required time to converge for each of the considered algorithms when applied to a network of 50 sensors, depicted in Figure \ref{fig:network}, with four different measurement noise standard deviation, namely $0.01, 0.05$, $0.1$ and~$0.3$. \normalsize}
\label{fig:ex13}
\vspace{-5mm}
\end{center}
\end{figure}
The computational time for the considered algorithms are presented in Figure \ref{fig:ex13}. As can be seen from the figure DDRA is at least twice as fast as DPDLA, owing to very simple computations required from each agent at every iteration. This is the case if both algorithms are executed in a centralized manner and if we neglect the communication cost or delay. Based on the presented results, our proposed algorithm provides more accurate estimates, and even though slower when implemented in a centralized manner, it provides a better distributed algorithm as it requires far less amount of communications.
\begin{figure}[t]
\begin{center}
\includegraphics[width=7.5cm]{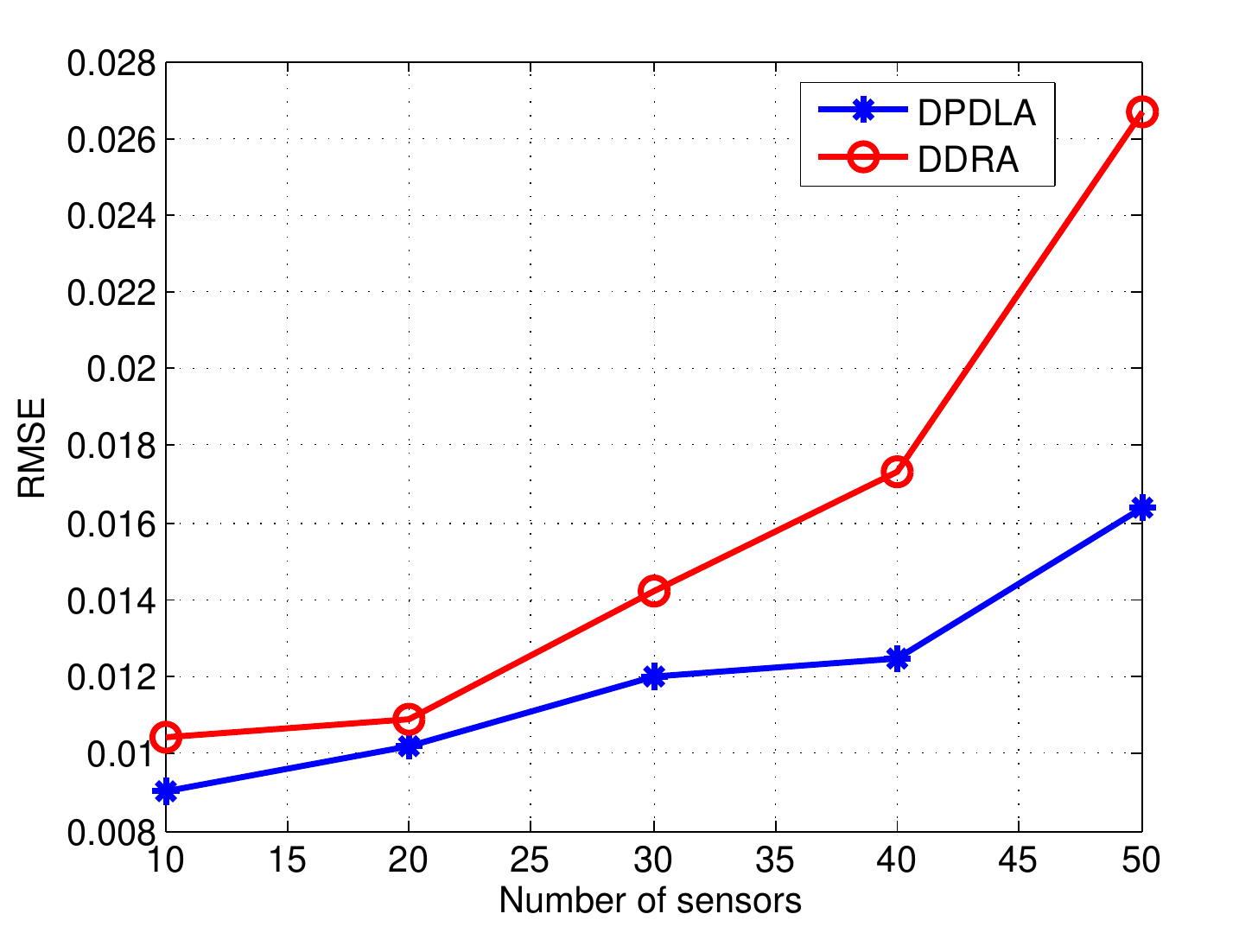}    
\caption{\footnotesize The RMSE results from the considered algorithms when applied to networks of varying number of sensors, with measurement noise standard deviation of 0.01.\normalsize}
\label{fig:ex21}
\vspace{-7mm}
\end{center}
\end{figure}

In the second simulation setup, we test the performance of the considered algorithms, when applied to networks with varying number of sensors, namely, 10, 20, 30, 40 and 50. In this setup we assume that the measurement noise standard deviation is 0.01, and we consider 50 instances for each network size. Furthermore the size of the considered area and the communication range, $r_c$, for each network size are chosen such that the resulting inter-sensor measurement graphs are connected but loosely. Figure \ref{fig:ex21} illustrates the RMSE results for this experiment. As before, as can be seen from the figure, DPDLA provides more accurate estimates for all network sizes.
\begin{figure}[t]
\begin{center}
\includegraphics[width=7.5cm]{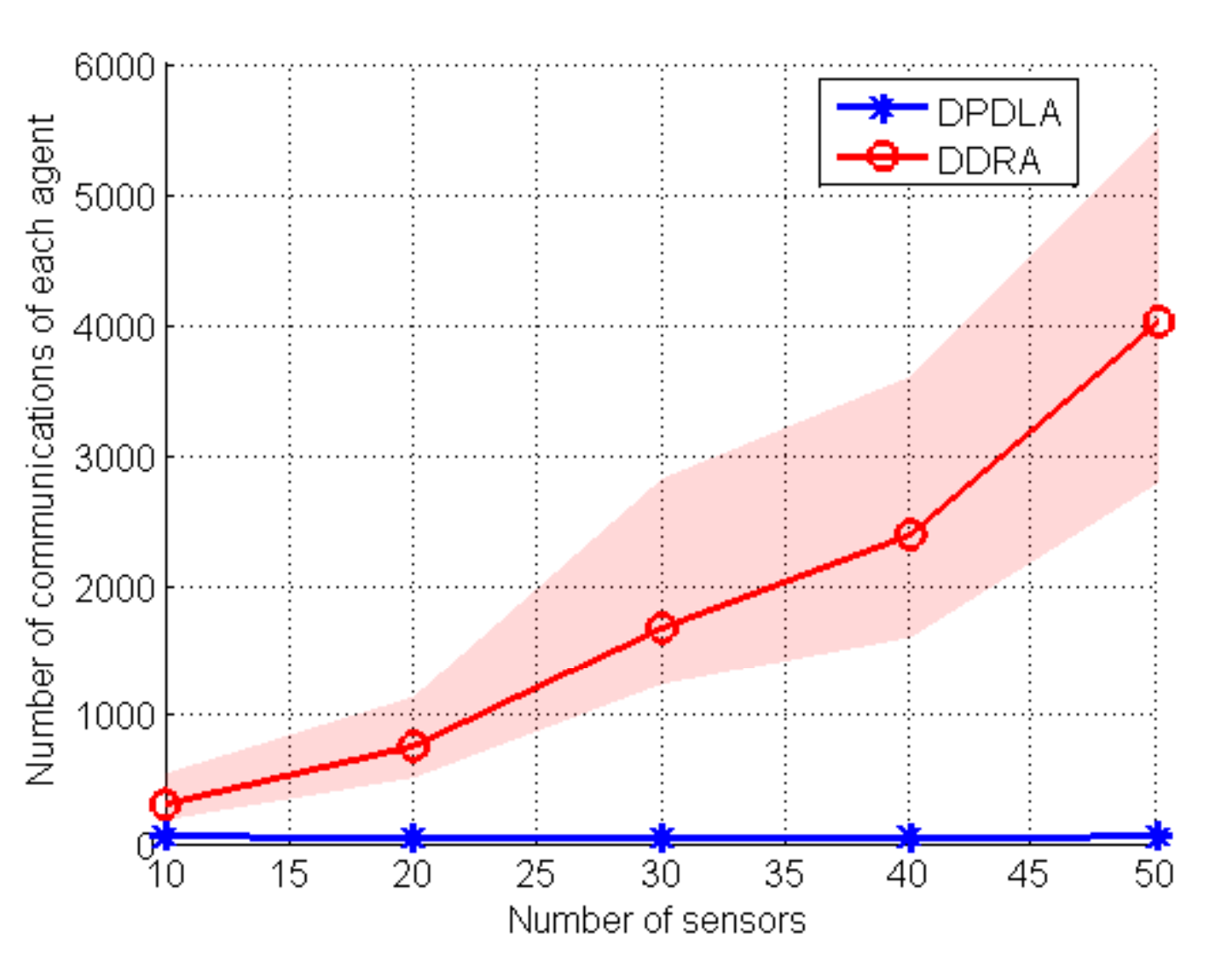}    
\caption{\footnotesize The number of communications that each agent needs to conduct for each of the algorithms to converge to a solution when applied to networks of varying number of sensors, with measurement noise standard deviation of 0.01. \normalsize}
\label{fig:ex22}
\vspace{-8mm}
\end{center}
\end{figure}
\begin{figure}[t]
\begin{center}
\includegraphics[width=7.5cm]{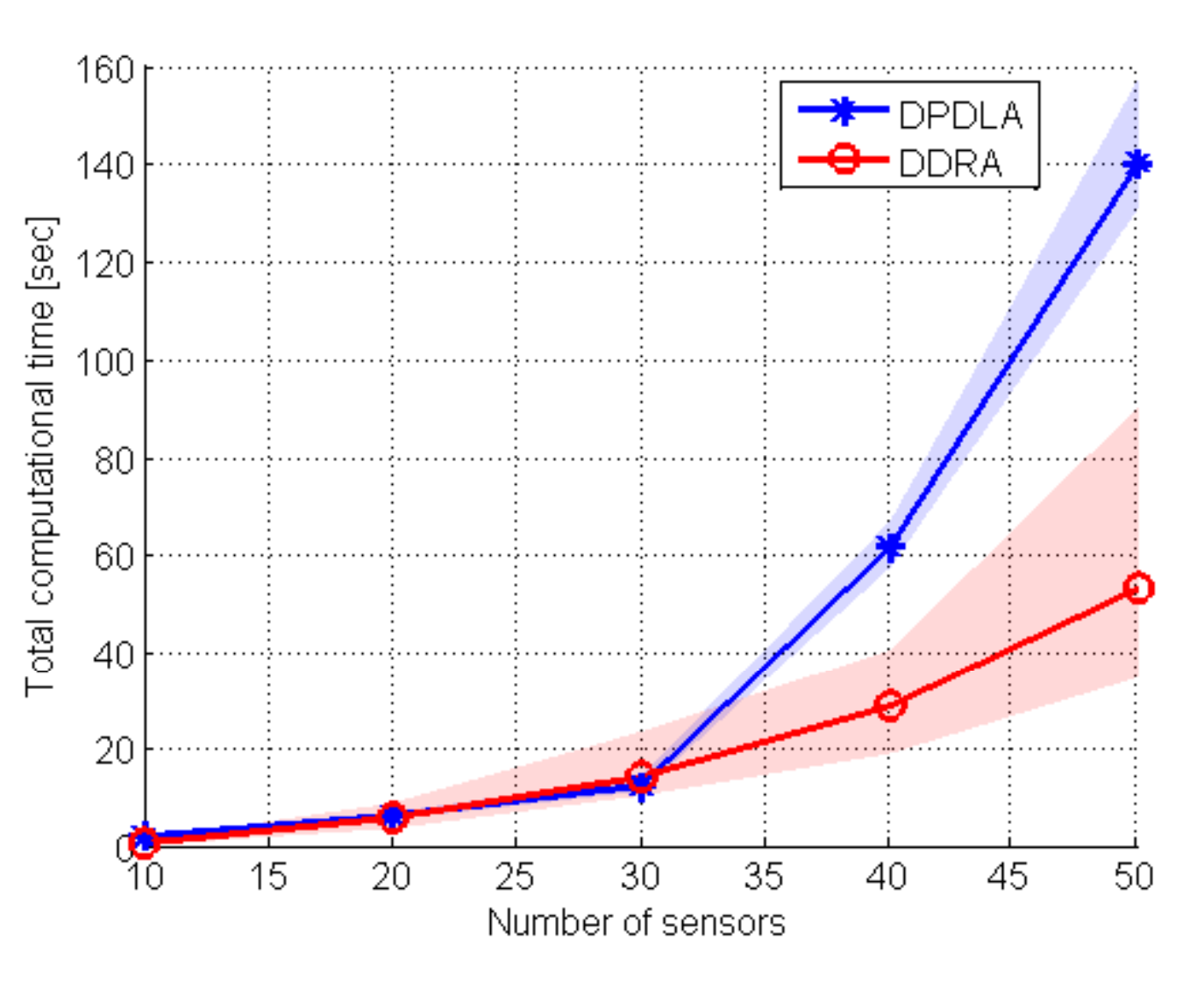}    
\caption{\footnotesize The required time for each of the considered algorithms to converge when applied to networks of varying number of sensors, with measurement noise standard deviation of 0.01. \normalsize}
\label{fig:ex23}
\vspace{-7mm}
\end{center}
\end{figure}
Also as can be seen from Figure \ref{fig:ex22}, the estimates are computed using far fewer communications among agents. The primal-dual method converged within around 11 iterations and the heights of the clique trees for the different sensor networks ware between 3 to 8. As can be seen from the figure, the number of required communications for the DDRA to converge grows much faster with network size than that of DPDLA which seems to be far less sensitive to this change. Figure \ref{fig:ex23} illustrates the total computational time of both algorithms when implemented in a centralized manner. As can be seen from this figure, our proposed algorithm requires similar or less amount of time to converge to a solution for networks of up to 30 sensors. Consequently, for networks with less than 30 sensors, our proposed algorithm outperforms DDRA in all the performance criteria. It is also worth mentioning that, the performance of our algorithm can be improved considerably, if the clustering of the sensors and generation of a clique tree are done using more sophisticated and tailored approaches. However, since we did not discuss such approaches, we abstained from any manipulation of the cliques and the clique tree and simply relied on standard and simple heuristics for this purpose, see e.g., \cite{kho:15c} and references~therein.

\subsection{Experiments Using Real Data}
\begin{figure}[t]
\begin{center}
\includegraphics[width=7.5cm]{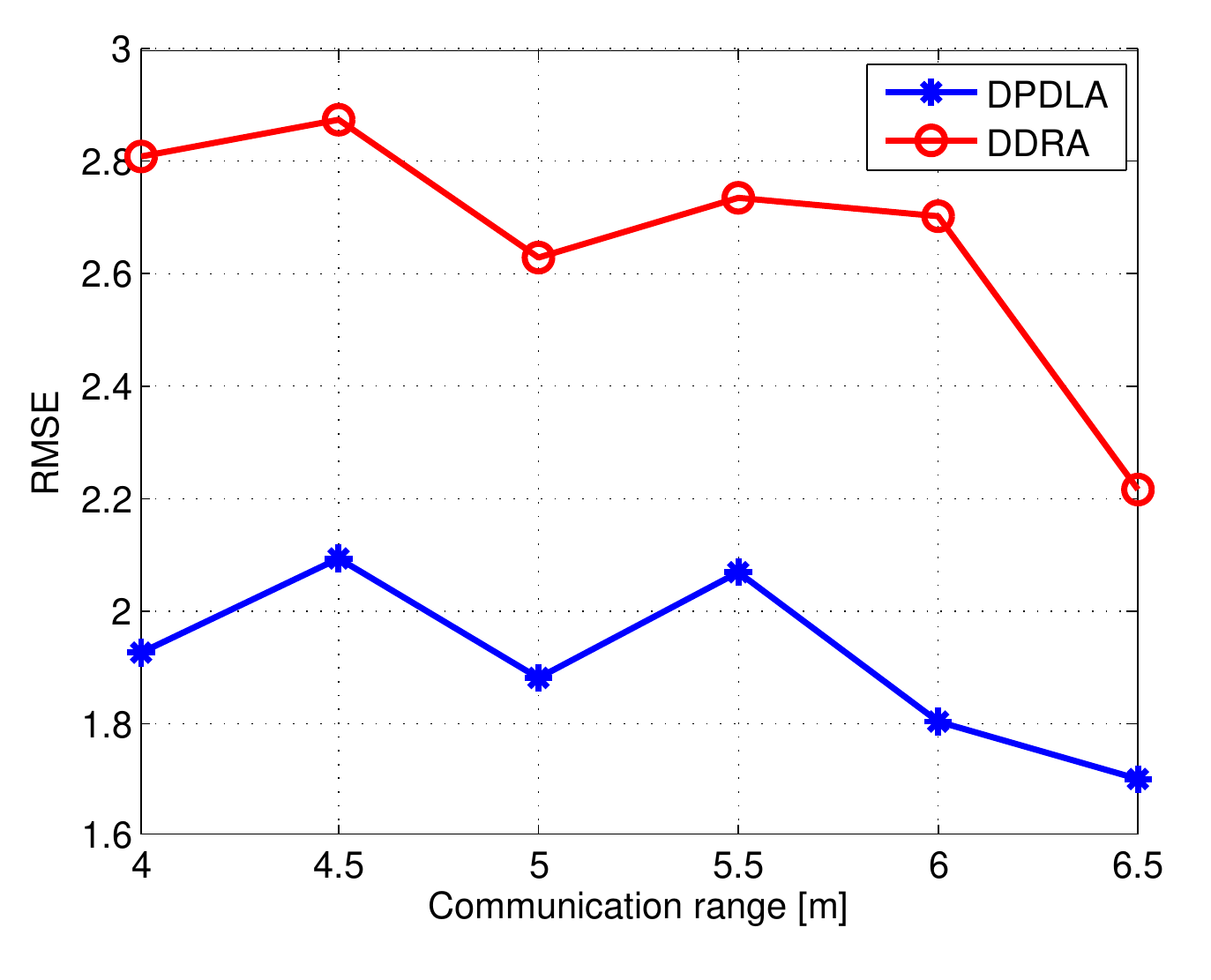}    
\caption{\footnotesize The RMSE results from the considered algorithms when applied to a localization problem based on real data, with a varying communication range. The $*$-marked line illustrates the RMSE results from DPDLA, whereas the o-marked line shows the RMSE results from DDRA.\normalsize}
\label{fig:ex31}
\vspace{-3mm}
\end{center}
\end{figure}

\begin{figure}[t]
\begin{center}
\includegraphics[width=7.5cm]{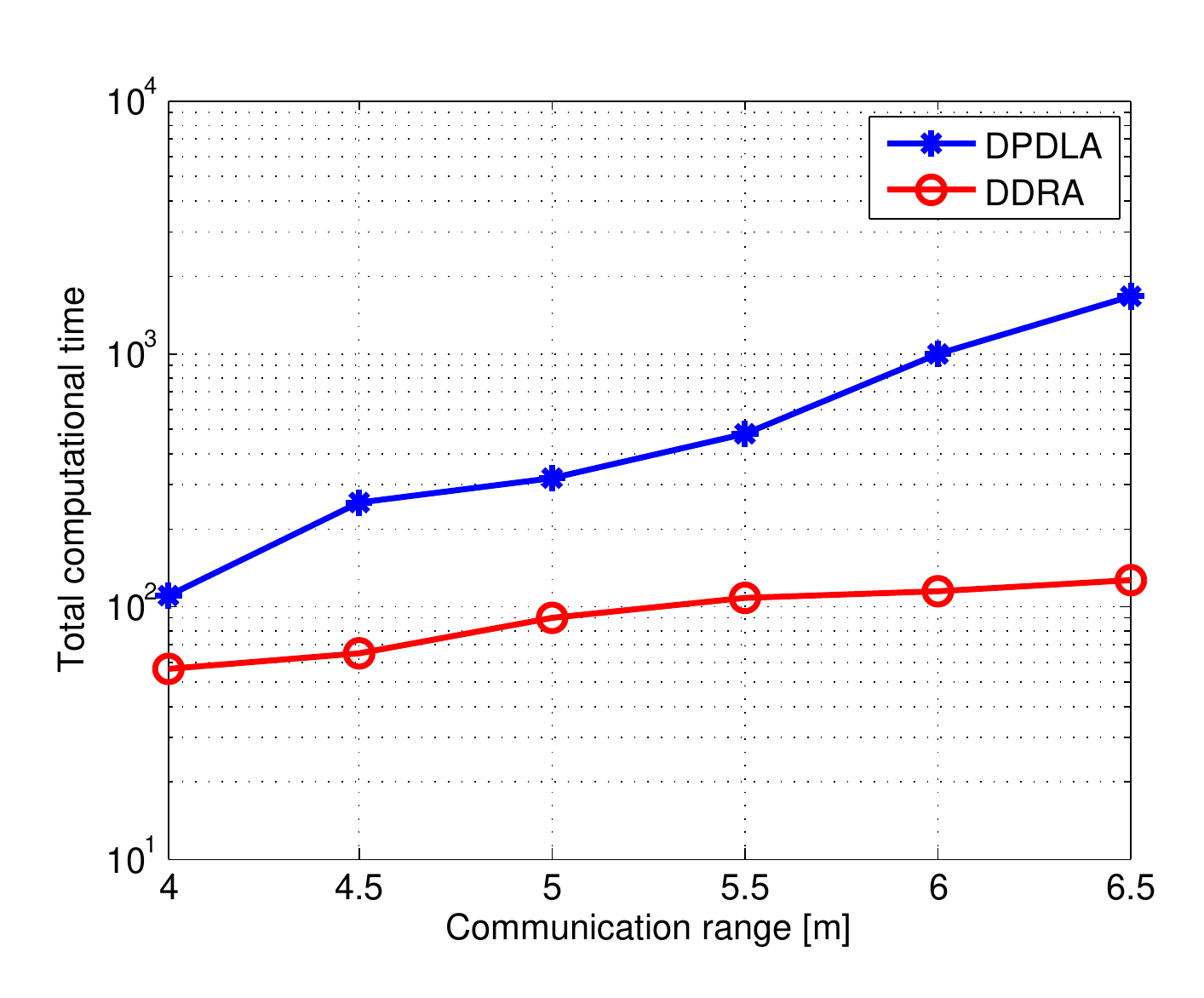}    
\caption{\footnotesize The required time for each of the considered algorithms to converge when applied to a localization problem based on real data, with a varying communication range.\normalsize}
\label{fig:ex33}
\vspace{-5mm}
\end{center}
\end{figure}

In this section, we present the results from conducted experiments based on real data. This data was taken from \cite{pat:03}, that includes time of arrival (TOA) measurements among 44 sensors, 4 of which are deemed to be anchors. The sensors are spread out in a $14 \times 13$ area. We extract the range measurements from the available TOA measurements. This provides us with biased range measurements with a standard deviation of 1.82 meters, see \cite{pat:03}. We here study the performance of DDRA and DPDRA for different levels of connectivity of the inter-sensor graph. To this end, we gradually change the communication range from 4 to 6.5 meters. Figures \ref{fig:ex31} and \ref{fig:ex33} illustrate the results.
Notice that due to biasedness and quality of the measurements, the intersection of the range measurement disks can be empty and hence DDRA fails to converge. This is because the gradient of the cost function of the disk relaxation problem does not vanish. Consequently, this algorithm has been terminated after 5000 iterations. Figure \ref{fig:ex31} illustrates the RMSE results from the experiment, which clearly depicts that DPDLA outperforms DDRA.
Furthermore, DPDLA required each agent was required to communicate with its neighbors around 100 times which seemed to be robust with respect to the level of connectivity of the inter-sensor range measurement graph. The primal-dual method for all these instances converged within roughly 17 iterations and the height of the clique tree varied between 7 to 9. Figure~\ref{fig:ex33} illustrates the computational time for DDRA and DPDLA. As was also observed from the experiments in Section \ref{sec:numerical-A}, DDRA clearly outperforms DPDLA when implemented in a centralized manner.
\section{Conclusions}\label{sec:conclusions}
In this paper we proposed a distributed localization algorithm for tree-structured scattered sensor networks founded on semidefinite relaxation of the localization problem. This algorithm is based on state-of-the-art primal-dual interior-point methods and relies on message-passing or dynamic programming over trees to distribute the computations. Due to this, the resulting algorithm requires far fewer steps and even fewer communications among computational agents to converge to an accurate solution, and it achieves this by putting a moderate computational burden on the agents. Furthermore, the proposed distributed algorithm is robust to biases in the measurements, or in general bad quality of the measurements. This stems from the power of semidefinite relaxation for localization problems. Despite these advantages, the proposed algorithm is much more complicated than algorithms that rely on first-order methods. This is largely due to the fact that generally second-order methods are far more complicated than their first-order counter parts.

The choice of clustering of the sensors and the strategy for assigning the available measurements to computational agents can have a significant effect on the performance of our proposed algorithm. Also smart clustering of the sensors, may even enable us to use the computational infrastructure at the anchors and utilize them  as computational agents. In this paper, we briefly discussed the importance of this and provided some suggestions on how the used heuristic strategies for this purpose can be improved. We did not investigate this topic in detail, however, we believe that further exploration of this matter can result in interesting results. Furthermore, distributed approaches for computing cliques and clique trees of the inter-sensor measurement were not covered in this paper, although, complementing the proposed algorithm with such methods can enhance the practicality of the algorithm.
\vspace{-10pt}
%

%
%

\ifCLASSOPTIONcaptionsoff
  \newpage
\fi



%
\bibliographystyle{plain}
\bibliography{IEEETrans}

\begin{thebibliography}{10}

\bibitem{and:11}
M.~S. Andersen.
\newblock {\em Chordal Sparsity in Interior-Point Methods for Conic
  Optimization}.
\newblock Ph{D} dissertation, university of {California, Los Angeles}, 2011.

\bibitem{ber:73}
U.~Bertelè and F.~Brioschi.
\newblock On non-serial dynamic programming.
\newblock {\em Journal of Combinatorial Theory, Series A}, 14(2):137--148,
  1973.

\bibitem{bis:06}
P.~Biswas, T.-C. Liang, K.-C. Toh, Y.~Ye, and T.-C. Wang.
\newblock Semidefinite programming approaches for sensor network localization
  with noisy distance measurements.
\newblock {\em IEEE Transactions on Automation Science and Engineering},
  3(4):360--371, Oct 2006.

\bibitem{bis:04}
P.~Biswas and Y.~Ye.
\newblock Semidefinite programming for ad hoc wireless sensor network
  localization.
\newblock In {\em Proceedings of the 3rd international symposium on Information
  processing in sensor networks}, pages 46--54. ACM, 2004.

\bibitem{blp:94}
J.~R.~S. Blair and B.~W. Peyton.
\newblock An introduction to chordal graphs and clique trees.
\newblock In J.~A. George, J.~R. Gilbert, and J.~W-H. Liu, editors, {\em Graph
  Theory and Sparse Matrix Computations}, volume~56, pages 1--27.
  Springer-Verlag, 1994.

\bibitem{bul:00}
N.~Bulusu, J.~Heidemann, and D.~Estrin.
\newblock {GPS}-less low-cost outdoor localization for very small devices.
\newblock {\em IEEE Personal Communications}, 7(5):28--34, Oct 2000.

\bibitem{cha:09}
F.~Chan and H.-C. So.
\newblock Accurate distributed range-based positioning algorithm for wireless
  sensor networks.
\newblock {\em IEEE Transactions on Signal Processing}, 57(10):4100--4105,
  2009.

\bibitem{fukuda_exploitingsparsity}
M.~Fukuda, M.~Kojima, , K.~Murota, and K.~Nakata.
\newblock Exploiting sparsity in semidefinite programming via matrix completion
  {I}: General framework.
\newblock {\em SIAM Journal on Optimization}, 11:647--674, 2000.

\bibitem{gho:13}
M.~R. Gholami, L.~Tetruashvili, E.~G. Strom, and Y.~Censor.
\newblock Cooperative wireless sensor network positioning via implicit convex
  feasibility.
\newblock {\em IEEE Transactions on Signal Processing}, 61(23):5830--5840,
  2013.

\bibitem{gol:04}
M.~C. Golumbic.
\newblock {\em Algorithmic Graph Theory and Perfect Graphs}.
\newblock Elsevier, 2nd edition, 2004.

\bibitem{gro:84}
R.~Grone, C.~R. Johnson, E.~M. \'{S}, and H.~Wolkowicz.
\newblock Positive definite completions of partial hermitian matrices.
\newblock {\em Linear Algebra and its Applications}, 58:109--124, 1984.

\bibitem{kho:15c}
S.~Khoshfetrat~Pakazad, A.~Hansson, and M.~S. Andersen.
\newblock Distributed primal-dual interior-point methods for solving
  tree-structured coupled problems using message passing.
\newblock {\em Optimization Methods and Software}, July 2016.

\bibitem{kim+koj+mev+yam10}
S.~Kim, M.~Kojima, M.~Mevissen, and M.~Yamashita.
\newblock Exploiting sparsity in linear and nonlinear matrix inequalities via
  positive semidefinite matrix completion.
\newblock {\em Mathematical Programming}, pages 1--36, 2010.

\bibitem{kim:09}
S.~Kim, M.~Kojima, and H.~Waki.
\newblock Exploiting sparsity in {SDP} relaxation for sensor network
  localization.
\newblock {\em SIAM Journal on Optimization}, 20(1):192--215, 2009.

\bibitem{kol:09}
D.~Koller and N.~Friedman.
\newblock {\em Probabilistic Graphical Models: Principles and Techniques}.
\newblock MIT press, 2009.

\bibitem{mor:97}
J.~J. Moré and Z.~Wu.
\newblock Global continuation for distance geometry problems.
\newblock {\em SIAM Journal on Optimization}, 7(3):814--836, 1997.

\bibitem{nar:14}
M.~Naraghi-Pour and G.~Rojas.
\newblock A novel algorithm for distributed localization in wireless sensor
  networks.
\newblock {\em ACM Transactions on Sensor Networks}, 11(1):1, 2014.

\bibitem{kho:15d}
S.~Khoshfetrat Pakazad, A.~Hansson, M.~S. Andersen, and A.~Rantzer.
\newblock Distributed semidefinite programming with application to large-scale
  system analysis.
\newblock {\em ArXiv e-prints}, April 2015.

\bibitem{pat:03}
N.~Patwari, A.~O. Hero~III, M.~Perkins, N.~Correal, and R.~J. O'dea.
\newblock Relative location estimation in wireless sensor networks.
\newblock {\em IEEE Transactions on Signal Processing}, 51(8):2137--2148, 2003.

\bibitem{rec:10}
B.~Recht, M.~Fazel, and P.~A. Parrilo.
\newblock Guaranteed minimum-rank solutions of linear matrix equations via
  nuclear norm minimization.
\newblock {\em SIAM Review}, 52(3):471--501, 2010.

\bibitem{sch:15}
S.~Schlupkothen, G.~Dartmann, and G.~Ascheid.
\newblock A novel low-complexity numerical localization method for dynamic
  wireless sensor networks.
\newblock {\em IEEE Transactions on Signal Processing}, 63(15):4102--4114, Aug
  2015.

\bibitem{shi:10}
Q.~Shi, C.~He, H.~Chen, and L.~Jiang.
\newblock Distributed wireless sensor network localization via sequential
  greedy optimization algorithm.
\newblock {\em IEEE Transactions on Signal Processing}, 58(6):3328--3340, June
  2010.

\bibitem{sim:14}
A.~Simonetto and G.~Leus.
\newblock Distributed maximum likelihood sensor network localization.
\newblock {\em IEEE Transactions on Signal Processing}, 62(6):1424--1437, March
  2014.

\bibitem{soa:15}
C.~Soares, J.~Xavier, and J.~Gomes.
\newblock Simple and fast convex relaxation method for cooperative localization
  in sensor networks using range measurements.
\newblock {\em IEEE Transactions on Signal Processing}, 63(17):4532--4543, Sep
  2015.

\bibitem{sri:08}
S.~Srirangarajan, A.~H Tewfik, and Z.~Luo.
\newblock Distributed sensor network localization using {SOCP} relaxation.
\newblock {\em Wireless Communications, IEEE Transactions on},
  7(12):4886--4895, 2008.

\bibitem{tod:96}
M.~J. Todd, K.~C. Toh, and R.~H. T{\"u}t{\"u}nc{\"u}.
\newblock On the {Nesterov-Todd} direction in semidefinite programming.
\newblock {\em SIAM Journal on Optimization}, 8:769--796, 1996.

\bibitem{wan:06}
Z.~Wang, S.~Zheng, S.~Boyd, and Y.~Ye.
\newblock Further relaxations of the sdp approach to sensor network
  localization.
\newblock Technical report, Stanford University, Tech. Rep, 2006.

\bibitem{wri:97}
S.~J. Wright.
\newblock {\em Primal-Dual Interior-Point Methods}.
\newblock Society for Industrial and Applied Mathematics, 1997.

\end{thebibliography}

%





\end{document}